\documentclass[aos,MSNbibl,dvips]{arximspdf}
\usepackage{ifthen}
\usepackage{mathbh}
\usepackage{twoopt}

\doi{10.1214/12-AOS1047} 
\volume{40}
\issue{5}
\pubyear{2012}
\firstpage{2697}
\lastpage{2732}

\makeatletter

\newcommand{\rright}{\right}
\newcommand{\lleft}{\left}
\newcommand{\rrVert}{\Vert}
\newcommand{\rrvert}{\vert}
\newcommand{\llVert}{\Vert}
\newcommand{\llvert}{\vert}

\newtheorem{theorem}{Theorem}
\newtheorem{prop}[theorem]{Proposition}
\newtheorem{cor}[theorem]{Corollary}
\newtheorem{lem}[theorem]{Lemma}

\newproclaim{defi}{Definition}
\newproclaim{rem}{Remark}

\newcommand{\iint}{\int\!\!\int}

\newcommand{\nset}{\mathbb{N}}
\newcommand{\rset}{\mathbb{R}}
\newcommand{\zset}{\mathbb{Z}}

\def\rset{\mathbb{R}}
\def\Xset{\mathsf{X}}
\def\Xsetprod{\bar{\Xset}}
\def\Xsigma{\mathcal{X}}
\def\Xsigmaprod{\bar{\Xsigma}}

\newcommand{\Yset}{\mathsf{Y}}
\newcommand{\Ysigma}{\mathcal{Y}}
\newcommand{\Pprod}{\bar{P}}

\def\Zset{\mathsf{Z}}

\def\Kset{\mathsf{K}}

\def\Mset{\mathsf{M}}
\def\Oset{\mathsf{O}}
\def\Cset{{\mathsf{C}}}
\def\Csetprod{{\bar{\Cset}}}

\def\Aset{\mathsf{A}}
\def\Asetprod{\bar{\Aset}}
\def\Dset{\mathsf{D}}

\newcommand{\one}{\mathbh{1}}

\newcommand{\PE}[2][]{
\ifthenelse{\equal{#1}{}}
{\ifthenelse{\equal{#2}{}}{\mathbb{E}}{\mathbb{E}[#2]}}
{\ifthenelse{\equal{#2}{}}{\mathbb{E}_{#1}}{\mathbb{E}_{#1}[#2]}}
}

\newcommand{\pblock}[2][]
{\ifthenelse{\equal{#1}{}}{\bolds{\ell}\langle#2\rangle}{\bolds{\ell}^{#1}\langle #2\rangle}
}

\newcommand{\PP}[2][]{
\ifthenelse{\equal{#1}{}}
{\ifthenelse{\equal{#2}{}}{\mathbb{P}}{\mathbb{P}[#2]}}
{\ifthenelse{\equal{#2}{}}{\mathbb{P}_{#1}}{\mathbb{P}_{#1}[#2]}}
}

\newcommand{\Pblockprod}[2][]
{\ifthenelse{\equal{#1}{}}{\bar{\mathbf{L}}\langle #2\rangle}{\bar{\mathbf{L}}^{#1} \langle #2 \rangle}
}

\def\rmd{\mathrm{d}}

\newcommand{\chunk}[4][]{\ifthenelse{\equal{#1}{}}{{{#2}_{#3}^{#4}}}{{#2}_{#1,#3}^{#4}}}

\newcommand{\Q}[1][]{\ifthenelse{\equal{#1}{}}{Q}{Q^{#1}}}

\newcommand{\q}[1][]{\ifthenelse{\equal{#1}{}}{q}{q^{#1}}}

\newcommand{\g}[1][]{\ifthenelse{\equal{#1}{}}{g}{g^{#1}}}

\newcommandtwoopt{\dens}[3][][]{
\ifthenelse{\equal{#1}{}}
{\ifthenelse{\equal{#2}{}}{p(#3)}{p_{#2}(#3)}}
{\ifthenelse{\equal{#2}{}}{p_{#1}(#3)}{p_{#1}^{#2}(#3)}}
}

\newcommandtwoopt{\cdens}[4][][]%
{
\ifthenelse{\equal{#1}{}}
{\ifthenelse{\equal{#2}{}}{p(#3|#4)}{p_{#2}(#3|#4)}}
{\ifthenelse{\equal{#2}{}}{p_{#1}(#3|#4)}{p_{#1}^{#2}(#3|#4)}}
}

\newcommand{\xprod}{\bar{x}}

\newcommand{\wprod}{\bar{w}}

\newcommand{\fprod}{\bar{f}}

\newcommand{\Pblock}[2][]
{\ifthenelse{\equal{#1}{}}{\mathbf{L}\langle#2\rangle}{\mathbf{L}^{#1}\langle#2\rangle}}

\newcommand{\DDelta}[4][]
{\ifthenelse{\equal{#1}{}}{\Delta_{#2} \langle #3 \rangle(#4)}{\Delta^{#1}_{#2} \langle #3 \rangle(#4)}}

\newcommand{\mle}[2]{\hat\theta_{#1,#2}}
\newcommand{\argmax}{\mathop{\arg\max}}

\newcommand{\URoot}[1][]{\ifthenelse{\equal{#1}{}}{R}{R_{#1}}}

\newcommand{\VRoot}[1][]{\ifthenelse{\equal{#1}{}}{S}{S_{#1}}}

\newcommand{\UCov}[1][]{\ifthenelse{\equal{#1}{}}{\URoot{^t\URoot}}{\URoot[#1] {^t\URoot[#1]}}}

\newcommand{\VCov}[1][]%
{%
\ifthenelse{\equal{#1}{}}{\VRoot{^t\VRoot}}{\VRoot[#1] {^t\VRoot[#1]}}%
}

\newcommand{\ocint}[1]{(#1]}
\newcommand{\ooint}[1]{(#1)}
\newcommand{\ccint}[1]{[#1]}

\newcommand{\prob}[3][]
{
\ifthenelse{\equal{#1}{}}{\lambda_{#2}\langle #3\rangle}{\lambda^{#1}_{#2}\langle #3\rangle}
}

\newcommand{\probprod}[3][]
{
\ifthenelse{\equal{#1}{}}{\bar{\lambda}_{#2} \langle3 \rangle}{\bar{\lambda}^{#1}_{#2} \langle #3 \rangle}
}

\makeatother

\begin{document}
\begin{frontmatter}

\title{Asymptotic properties of the maximum likelihood estimation in
misspecified hidden Markov models\thanksref{BOTH}}
\runtitle{MLE in misspecified HMMs}

\thankstext{BOTH}{Supported by the Agence Nationale de la Recherche
through the 2009-2012 project Big MC.}

\begin{aug}
\author[A]{\fnms{Randal} \snm{Douc}\corref{}\ead[label=e1]{randal.douc@it-sudparis.eu}}
\and
\author[B]{\fnms{Eric} \snm{Moulines}\ead[label=e2]{eric.moulines@telecom-paristech.fr}}
\runauthor{R. Douc and E. Moulines}
\affiliation{T\'el\'ecom SudParis and T\'el\'ecom ParisTech}
\address[A]{SAMOVAR\\
CNRS UMR 5157\\
Institut T\'el\'ecom/T\'el\'ecom SudParis\\
9 rue Charles Fourier\\
91000 Evry\\
France\\
\printead{e1}} 
\address[B]{LTCI\\
CNRS UMR 5141\\
Institut T\'el\'ecom /T\'el\'ecom ParisTech\\
46 rue Barrault\\
75634 Paris Cedex 13\\
France\\
\printead{e2}}
\end{aug}

\received{\smonth{10} \syear{2011}}
\revised{\smonth{7} \syear{2012}}

%
\begin{abstract}
Let $(Y_k)_{k \in\zset}$ be a stationary sequence on a probability
space $(\Omega,\mathcal{A},\PP{})$ taking values in a standard Borel
space $\Yset$. Consider the associated maximum likelihood estimator
with respect to a parametrized family of hidden Markov models such that
the law of the observations $(Y_k)_{k \in\zset}$ is \textit{not}
assumed to be described by any of the hidden Markov models of this
family. In this paper we investigate the consistency of this estimator
in such misspecified models under mild assumptions.
\end{abstract}

%

\end{frontmatter}

\section{Introduction}

An assumption underlying most of the classical theory of maximum
likelihood is that the ``true'' distribution of the observations is
known to lie within a specified parametric family of distributions. In
many settings, it is doubtful that this assumption is satisfied. It is
therefore natural to investigate the convergence of the maximum
likelihood estimator (MLE) and to identify the possible limit for
misspecified models. Such questions have been mainly investigated for
models in which observations are independent; see
\cite{huber1967,white1982}. Much less is known on the behavior of the
MLE estimate for dependent observations; see~\cite{fombyhill2003} and
the references therein.

For independent observations, under mild additional technical
conditions, the MLE converges to the parameter which minimizes the
relative entropy rate; see~\cite{huber1967}. The purpose of this paper
is to show that such a result remains true when the observations are
from an ergodic process and for classes of parametric distributions
associated to hidden Markov models (HMM). A HMM is a bivariate
stochastic process $(X_k, Y_k)_{k\geq0}$, where $(X_k)_{k\geq0}$ is a
Markov chain (often referred to as the state sequence) in a state space
$\Xset$ and, conditionally on $(X_k)_{k \geq0}$, $(Y_k)_{k\geq0}$ is a
sequence of independent random variables in a state space $\Yset$ such
that the conditional distribution of $Y_k$ given the state sequence
depends on $X_k$ only. The key feature of HMMs is that the state
sequence $(X_k)_{k\geq0}$ is not observable, so that statistical
inference has to be carried out by means of the observations
$(Y_k)_{k\geq0}$ only. Such problems are far from straightforward due
to the fact that the observation process $(Y_k)_{k\geq0}$ is generally
a dependent, non-Markovian time series [despite that the bivariate
process $(X_k,Y_k)_{k\geq0}$ is itself a Markov chain].

HMMs have been intensively used in many scientific disciplines
including econometrics
\cite{hullwhite1987,mamonelliott2007}, biology
\cite{churchill1992}, engineering~\cite{juangrabiner1991},
neurophysiology~\cite{fredkinrice1987} and the statistical
inference is therefore of significant practical
importance~\cite{cappemoulinesryden2005}. In all these applications,
misspecified models are the rule,
so it is worthwhile to understand the behavior of MLE under such regime.

This work extends previous results in this direction obtained by Mevel
and Finesso
\cite{mevelfinesso2004}, but which are restricted to discrete
state-space Markov chains.
Our main result of consistency of the MLE in misspecified HMMs is
derived under assumptions which are quite weak,
covering general state-space HMMs under conditions which are much
weaker than~\cite{doucmoulinesryden2004}, where a strong mixing
condition was imposed on the transition kernels of the hidden chain.
Therefore our results can be applied to many models of practical
interest, including the Gaussian linear state space model, the discrete
state-space HMM and more general nonlinear state-space models.

The paper is organized as follows. In Section~\ref{secproblem-statement},
we first introduce the setting and notations that are used throughout
the paper.
In Section~\ref{secassumptions-mainresults}, we state our main
assumptions and results. In Section~\ref{secapplications}, our main
result is used to establish consistency in three general classes of
models: linear-Gaussian state space models, finite state models and
nonlinear state space models of the vector ARCH type (this includes the
stochastic volatility model and many other models of interest in time
series analysis and financial econometrics). Section~\ref{secproof}
is devoted to the proof of our main result.

\subsection*{Notation}
Some notation pertaining to transition kernels is required.
Let $L$ be a (possibly unnormalized) transition kernel on $(\Xset
,\Xsigma)$, that is, for any $x \in\Xset$, $L(x,\cdot)$ is a finite
measure on $(\Xset,\Xsigma)$ and for any $A \in\Xsigma$, $x \mapsto
L(x,A)$ is measurable function from $(\Xset,\Xsigma)$ to $
(\ccint
{0,1},\mathcal{B}(\ccint{0,1}) )$. $L$ acts on bounded
functions $f$ on $\Xset$ and on $\sigma$-finite positive measures
$\mu$ on
$(\Xset,\Xsigma)$ via
\[
L f(x) = \delta_x L f \triangleq\int L(x,\rmd y) f(y),\qquad \mu L (A) =
\mu L \one_A \triangleq\int\mu(\rmd x) L(x,A).
\]
If $L_1$ and $L_2$ are two transition kernels on $(\Xset,\Xsigma)$,
then $L_1 L_2$ is the transition kernel
on $(\Xset,\Xsigma)$, given, for any $x \in\Xset$ and $A \in
\Xsigma$ by
\[
L_1 L_2 (x,A)= \int L_1(x,\rmd y)
L_2(y, A).
\]
%
\section{Problem statement}
\label{secproblem-statement}

We consider a parameterized family of HMMs with
parameter space $\Theta$, assumed to be a compact metric space.
For each parameter $\theta\in\Theta$, the distribution of the HMM
is specified by the transition kernel $\Q[\theta]$ of
the Markov chain $(X_k)_{k\ge0}$, and by the conditional distribution
$\g[\theta]$ of the observation $Y_k$ given the hidden state $X_k$,
referred to as
the likelihood of the observation.

For any $m \leq n$ and any sequence $\{a_k \}_{k \in\zset}$,
denote $\chunk{a}{m}{n}\triangleq(a_m,\ldots,a_n)$, and for any
probability
measure $\chi$ on $(\Xset,\Xsigma)$,
define the likelihood of the observations by
\begin{eqnarray*}
p_{\chi}^{\theta}\bigl(y_{m}^{n}\bigr)
&\triangleq&\int\cdots\int\chi
(\rmd
x_m) \g[\theta](x_m, y_m) \prod
_{p = m + 1}^n \Q[\theta](x_{p - 1}, \rmd
x_p) \g[\theta](x_p, y_p),
\\
p_{\chi}^{\theta}\bigl(y_{p}^{n}|\chunk{y} {m} {p-1}\bigr) &\triangleq&
 p_{\chi}^{\theta}\bigl(y_{m}^{n}\bigr)/p_{\chi}^{\theta}\bigl(y_{m}^{p-1}\bigr),\qquad m< p
\leq n,
\end{eqnarray*}
with the standard convention $\prod_{p = m}^n a_p = 1$ if $m >
n$.\vspace*{1pt}

Let $(\Omega,{\mathcal F},\PP{})$ be a probability space, and let
$(Y_k)_{k \in\zset}$ be a \textit{stationary ergodic} stochastic process
taking value in $(\Yset,\Ysigma)$. We denote by $\PP[Y]{}$ the image
probability of $\PP{}$ by $(Y_k)_{k \in\zset}$ on the product space
$(\Yset^\zset,\Ysigma^{\otimes\zset})$, and $\PE[Y]{}$ the associated
expectation. We stress that the distribution $\PP[Y]{}$ may or may not
belong to the parametric family of distributions specified by the
transition kernels $ \{ (\Q[\theta],\g[\theta]), \theta\in
\Theta
\}$. If $\PP[Y]{}$ does not belong to $\mathcal{G}$, the model is
said to be misspecified.

If $\chi$ is a probability measure $(\Xset,\Xsigma)$, we define the
maximum likelihood estimator (MLE) associated to the initial
distribution $\chi$ by
%
\begin{equation}
\label{eqdef-MLE} \mle{\chi} {n} \triangleq\argmax_{\theta\in\Theta} \ln
 p_{\chi}^{\theta}\bigl(Y_{0}^{n-1}\bigr).
\end{equation}

The study of asymptotic properties of the MLE in HMMs was initiated in
the seminal work of Baum and Petrie~\cite{baumpetrie1966,petrie1969}
in the 1960s. In these papers, the model is assumed to be well
specified, and
the state space $\Xset$ and the
observation space $\Yset$ were both presumed to be finite sets. More
than two decades later, Leroux~\cite{leroux1992} proved consistency for
well-specified models in the case that $\Xset$ is a finite set, and
$\Yset$ is a general state
space. The consistency of the MLE in more general HMMs has subsequently
been investigated for well-specified models in a series of
contributions \cite
{leglandmevel2000a,leglandmevel2000,doucmatias2002,doucmoulinesryden2004,genon-catalotlaredo2006}
using different methods. A general
consistency result for HMMs has
been developed in~\cite{doucmoulinesolssonvanhandel2011}.

Though the consistency results above differ in the details of their
proofs, all proofs have a common thread which serves also as the
starting point for this paper. Denote by
$\dens[\chi][\theta]{\chunk{Y}{0}{n}}$ the likelihood of the
observations $Y_0^n$ for the
HMM with parameter $\theta\in\Theta$ and initial distribution $\chi$.
The first step of the proof aims to establish that for any
$\theta\in\Theta$, there is a constant $\ell(\theta)$ such that
\[
\lim_{n \to\infty} n^{-1} \log p_{\chi}^{\theta}\bigl(Y_{0}^{n-1}\bigr)
= \lim_{n \to\infty} n^{-1} \mathbb{E}\bigl[\log
p_{\chi}^{\theta}\bigl(Y_{0}^{n-1}\bigr)\bigr]
= \ell(\theta),\qquad \PP{}\mbox{-a.s.}
\]
Up to an additive constant, $\theta\mapsto\ell(\theta)$
is the negated relative entropy rate between the distribution of the
observations
and $\dens[\chi][\theta]{\cdot}$, respectively.
When the model is well-specified and $\theta= \theta_\star$ is the true
value of the parameter,
this convergence follows from the generalized
Shannon--Breiman--McMillan theorem~\cite{barron1985};
for misspecified models or for well-specified models with $\theta\ne
\theta_\star$
the existence of the limit is far from obvious.

The second step of the proof aims to prove that the maximizer of the likelihood
$\theta\mapsto n^{-1} \log\dens[\chi][\theta]{\chunk{Y}{0}{n}}$ converges
$\PP{}$-a.s. to the maximizer of $\theta\mapsto\ell(\theta)$,
that is,
to the minimizer of the relative entropy rate.
Together, these two steps show that the MLE is a natural estimator for
the parameters
which minimizes the relative entropy rate in the parametric family \mbox{$\{
(\Q[\theta], \g[\theta]), \theta\in\Theta\}$}.

Let us note that one could write the likelihood as
\[
n^{-1}\log p_{\chi}^{\theta}\bigl(Y_{0}^{n-1}\bigr) =
\frac{1}{n}\sum_{k=0}^{n-1} \log p_{\chi}^{\theta}\bigl(Y_k| \chunk{Y} {0} {k-1}\bigr),
\]
where $\dens[\chi][\theta]{Y_k| \chunk{Y}{0}{k-1}}$ denotes the
conditional density of $Y_k$ given $Y_0^{k-1}$ under the misspecified
model with parameter $\theta$ (i.e., the one-step predictive density).
If the limit of $\dens[\chi][\theta]{Y_k| \chunk{Y}{0}{k-1}} \to
\pi^\theta_Y(\chunk{Y}{-\infty}{k})$ as $k\to\infty$ can be shown to
exist $\PP{}$-a.s. and to be $\PP{}$-integrable, the convergence of the
log-likelihood to the relative entropy rate follows from the Birkhoff
ergodic theorem, since the process $\{Y_k\}_{k \in\zset}$ is assumed
to be ergodic. This result provides an explicit representation of the
relative entropy rate $\ell(\theta)$ as the expectation of the limit
$\ell(\theta) = \PE{\log\pi^\theta_Y(\chunk{Y}{-\infty}{0})}$. The
limit $\pi^\theta_Y(\chunk{Y}{-\infty}{k})$ might be interpreted as the
conditional likelihood of $Y_k$ given the whole past
$\chunk{Y}{-\infty}{k-1}$, but we must refrain ourselves of considering
this quantity as a conditional density.

Such an approach was used in~\cite{baumpetrie1966} for finite state-space,
and was later extended by Douc, Moulines and Ryd\'{e}n \cite
{doucmoulinesryden2004} to general state-space, but under stringent
technical conditions
(uniform mixing of the Markov kernel, which more or less restricts the
validity of the results to compact state-spaces, leaving aside
important models, such as Linear Gaussian state-space models).

Alternatively, the
predictive distribution $\dens[\chi][\theta]{Y_k| \chunk{Y}{0}{k-1}}$
can be expressed
as a component of the state of a measure-valued Markov chain; in this approach,
the existence of the limiting relative entropy rate $\ell(\theta)$,
follows from
the ergodic theorem for Markov chains, provided that this Markov chain
can be shown to be ergodic.
This approach was used in \cite
{leglandmevel2000a,leglandmevel2000,doucmatias2002} and was later
extended to misspecified models by
White~\cite{mevelfinesso2004}. This technique
is adequate for finite state-space Markov chains, but does not extend
easily to general state-space Markov chains; see~\cite{doucmatias2002}.

In~\cite{leroux1992}, the existence of the relative entropy rate is
established by means of Kingman's subadditive ergodic theorem (the same
approach is used indirectly in~\cite{petrie1969}, which invokes the
Furstenberg--Kesten theory of random matrix products). After some
additional work, an explicit representation of the relative enropy rate
is again obtained. However, as is noted in~\cite{leroux1992}, page~136,
the latter is surprisingly difficult, as Kingman's ergodic theorem does
not directly yield a representation of the limit as an expectation.

For completeness, we note that a recent attempt~\cite{fuh2006} to prove
consistency of the MLE for general HMMs contains very serious problems in
the proof~\cite{jensen2009} (not addressed in~\cite{fuh2009}), and
therefore fails to establish the claimed results.

In this paper, we prove consistency of the MLE for general HMMs in
misspecified models under quite general
assumptions. Our proof follows broadly the original approach of Baum
and Petrie~\cite{baumpetrie1966} and Douc, Moulines and Ryd\'{e}n~\cite{doucmoulinesryden2004},
but relaxes the very restrictive
technical conditions used in these works and extends the analysis to
misspecified models.
The key technique to obtain this result is to establish the exponential
forgetting of the
filtering distribution; this result is obtained by using an original
coupling technique originally introduced in
\cite{kleptsynaveretennikov2008} and refined in \cite
{doucfortmoulinespriouret2009}.

\section{Assumptions and main results}
\label{secassumptions-mainresults} For any integer $t\geq1$, $\theta
\in\Theta$ and any sequence $\chunk{y}{0}{t-1} \in\Yset^t$, consider
the unnormalized kernel $\Pblock[\theta]{\chunk{y}{0}{t-1}}$ on
$(\Xset,\Xsigma)$ defined for all $x_0 \in\Xset$ and $\Aset\in
\Xsigma$, by
%
\begin{equation}
\label{eqdef-Pblock}
\mathbf{L}^{\theta}\bigl\langle\chunk{y} {0} {t-1} \bigr\rangle(x_0,
\Aset)= \int\cdots\int\Biggl[\prod_{i=0}^{t-1}
\g[\theta](x_{i},y_{i}) \Q[\theta](x_{i},
\rmd x_{i+1}) \Biggr] \one_\Aset(x_{t}).
\end{equation}
Note that, for any $t \geq1$, $\theta\in\Theta$, $x_0 \in\Xset$,
and $\chunk{y}{0}{t-1} \in\Yset^t$,
%
\begin{equation}
\label{eqP-theta-vraisemblance} \mathbf{L}^{\theta}\bigl\langle\chunk{y} {0}
{t-1} \bigr\rangle(x_0,\Xset)= p_{x_0}^{\theta}\bigl(y_{0}^{t-1}\bigr),
\end{equation}
where for $x \in\Xset$, $s\leq t$,
$\dens[x][\theta]{\chunk{y}{s}{t}}$, the likelihood of the observation
$\chunk{y}{s}{t}$ starting from state $x$, is a shorthand notation for
$\dens[\delta_x][\theta]{\chunk{y}{s}{t}}$.
%
\begin{defi} \label{defilocal-Doeblin}
Let $r$ be an integer. A set $\Cset\in\Xsigma$ is a $r$-\textit{local
Doeblin set} with respect to the family $\{\Q[\theta], \g[\theta]\}
_{\theta\in
\Theta}$, if there exist positive functions $\epsilon^-_\Cset\dvtx  \Yset^r
\!\to\!\rset^+$, $\epsilon^+_\Cset\dvtx  \Yset^r \!\to\!\rset^+$ and a family of
probability measures $\{
\prob[\theta]{\Cset}{z}\}_{\theta\in\Theta
, z
\in\Yset^r}$ and of positive functions $\{\varphi_{\Cset}^\theta
\langle z \rangle\}_{\theta\in\Theta,z \in\Yset^r}$ such that for
any $\theta\in\Theta$, $z \in\Yset^r$, $
\prob[\theta]{\Cset
}{z}(\Cset
) =1$ and, for any $\Aset\in\Xsigma$, and $x \in\Cset$,
%
\begin{equation}
\label{eqdefinition-LD-set}\quad \epsilon^{-}_\Cset(z)
\varphi_{\Cset}^\theta\langle z\rangle(x)
\prob[\theta]{\Cset}
{z} (\Aset) \leq\mathbf{L}^{\theta}\langle z \rangle(x, \Aset\cap\Cset) \leq\epsilon
^{+}_{\Cset}(z)
\varphi_{\Cset}^\theta\langle z\rangle(x)
\prob[\theta]{\Cset}
{z}(\Aset).
\end{equation}
\end{defi}
This implies that for any measurable nonnegative function $f$ on
$(\Xset
,\Xsigma)$, $x \in\Cset$ and any $z \in\Yset^r$,
\[
\epsilon^{-}_\Cset(z) \varphi_{\Cset}^\theta
\langle z \rangle(x)
\prob[\theta]{\Cset} {z} (\one_\Cset f) \leq
\delta_x \Pblock[\theta]{z} (\one_\Cset f) \leq
\epsilon^{+}_{\Cset}(z) \varphi_{\Cset
}^\theta
\langle z \rangle(x)
\prob[\theta]{\Cset} {z} (\one_\Cset f).
\]

We require that the condition is satisfied for any $\theta\in\Theta$,
but this is not a serious
restriction since $\Theta$ is assumed to be compact.
%
\begin{rem}
To illustrate this condition, consider the case $r=1$. Assume that
for some set $\Cset$, there exist positive constants $\epsilon^-_\Cset,
\epsilon^+_\Cset$ and a family of probability measures $\{
\lambda_{\Cset}^\theta\}_{\theta\in\Theta}$ such that for any $\theta\in
\Theta$, $\lambda_{\Cset}^\theta(\Cset) =1$ and, for any $\Aset
\in
\Xsigma$, and $x \in\Cset$,
\[
\epsilon^{-}_\Cset\lambda_{\Cset}^\theta(
\Aset) \leq\Q^\theta(x, \Aset\cap\Cset) \leq\epsilon^{+}_{\Cset}
\lambda_{\Cset}^\theta(\Aset) .
\]
Then,\vspace*{1pt} clearly $\Pblock[\theta]{y}(x,\Aset)=\g^\theta(x,y) \Q^\theta
(x,\Aset)$ satisfies (\ref{eqdefinition-LD-set}) where $\epsilon_\Cset
^-$ and $\epsilon_\Cset^+$ are positive \textit{constants}. In this case
$\Cset$ is a $1$-local Doeblin set with respect to $\Q^\theta$; see
\cite{doucfortmoulinespriouret2009} and \cite
{kleptsynaveretennikov2008}.
\end{rem}
%
\begin{rem}
Local Doeblin sets share some similarities with $1$-small set in the
theory of Markov chains over general
state spaces; see~\cite{meyntweedie1993}, Chapter~5. Recall that a set
$\Cset$ is $1$-small for the kernel $\Q[\theta]$, $\theta\in\Theta$
if there exists a probability measure
$\tilde{\lambda}^\theta_\Cset$ and a constant $\tilde{\epsilon
}_\Cset>
0$, such that $\tilde{\lambda}^\theta_{\Cset}(\Cset) =1$, and for all
$x \in\Cset$ and $A \in\Xsigma$, $\Q[\theta](x,A \cap\Cset)
\geq
\tilde{\epsilon}_\Cset\tilde{\lambda}^\theta_\Cset(A \cap\Cset
)$. In
particular, a local Doeblin set is $1$-small with $\tilde{\epsilon
}_\Cset= \epsilon^-_\Cset$ and $\tilde{\lambda}^\theta_\Cset=
\lambda^\theta_\Cset$.
The main difference stems from the fact that we impose both a lower and
an upper bound, and we impose that the minorizing and the majorizing
measures are the same.
\end{rem}

\renewcommand\theenumi{(A\arabic{enumi})}
\renewcommand\labelenumi{\theenumi}
\begin{enumerate}
\item\label{assumlikelihoodDrift}
There exist an integer $r\geq1$ and a set $\Kset\in\Ysigma^{\otimes r}$
such that:
\end{enumerate}

{\renewcommand\theenumi{(\roman{enumi})}
\renewcommand\labelenumi{\theenumi}
\begin{enumerate}[(A1)\hspace*{5pt}(iii)]
\item\label{itemcondition-L-K}
$ \PP{\chunk{Y}{0}{r-1} \in\Kset}>2/3$.
\item For all $\eta> 0$, there exists a $r$-local Doeblin set $\Cset
\in\Xsigma$
such that for all $\theta\in\Theta$ and for all $\chunk
{y}{0}{r-1}\in\Kset$,
%
\begin{equation}
\label{eqbound-eta-G} \sup_{x_0 \in\Cset^c} p_{x_0}^{\theta}\bigl(y_{0}^{r-1}\bigr)
\leq\eta\sup_{x_0 \in\Xset}  p_{x_0}^{\theta}\bigl(y_{0}^{r-1}\bigr)<\infty
\end{equation}
and
%
\begin{equation}
\label{eqlower-bound} \inf_{\chunk{y}{0}{r-1} \in\Kset} \frac{\epsilon
_\Cset^-(\chunk
{y}{0}{r-1})}{\epsilon_\Cset^+(\chunk{y}{0}{r-1})} > 0,
\end{equation}
where the functions $\epsilon^+_\Cset$ and $\epsilon^-_\Cset$ are
defined in Definition~\ref{defilocal-Doeblin}.
\item\label{itemcondition-minoration}
There exists a set $\Dset$ such that
%
\begin{equation}
\label{eqcondition-minoration} \mathbb E\Bigl[{\ln^-
\inf_{\theta\in\Theta} \inf_{x \in\Dset}
\mathbf{L}^{\theta}\bigl\langle \chunk{Y} {0} {r-1}
\bigr\rangle(x,\Dset)}\Bigr] < \infty.
\end{equation}
\end{enumerate}}
%

\begin{enumerate}
\setcounter{enumi}{1}
\item\label{assummajo-g}
\end{enumerate}

{\renewcommand\theenumi{(\roman{enumi})}
\renewcommand\labelenumi{\theenumi}
\begin{enumerate}[(A2)\hspace*{5pt}(iii)]\vspace*{-26pt}
\item\label{itemassum-pos-g} For any $\theta\in\Theta$, the
function $\g[\theta]\dvtx  (x,y) \in\Xset\times\Yset\mapsto\g[\theta
](x,y)$ is positive,
\item\label{itemassum-supg} $\PE{ \ln^+\sup_{\theta\in\Theta}
\sup_{x \in\Xset} \g[\theta](x,Y_0)}<\infty$.
\end{enumerate}}

\begin{enumerate}
\setcounter{enumi}{2}
\item\label{assumcontinuity}
There exists $p \in\nset$ such that for any $x \in\Xset$ and $n
\geq
p$, $\PP{}$-a.s. the function $\theta\mapsto\dens[x][\theta
]{\chunk{Y}{0}{n}}$
is continuous on $\Theta$.
\end{enumerate}
%
%
\begin{rem}
Assumption~\ref{assummajo-g} assumes that the conditional likelihood $\g[\theta]$
is positive. The case where $\g[\theta]$ can vanish typically
requires
different conditions; see~\cite{budhirajaocone1997,vanhandel2008}. The second condition can be read as a generalized
moment condition on~$Y_0$. It is satisfied in many examples of interest.\vadjust{\goodbreak}
\end{rem}
%
\begin{rem}
To check~\ref{assumlikelihoodDrift}\ref{itemcondition-minoration},
one may, for example, check that:

\renewcommand\thelonglist{(\roman{longlist})}
\renewcommand\labellonglist{\thelonglist}
\begin{longlist}
\item\label{itemassummino-gQ} $\inf_{x \in\Dset} \inf_{\theta
\in
\Theta} \Q[\theta](x,\Dset) > 0$;
\item\label{itemassummino-gsuppDenom} $\PE{ \ln^- \inf_{\theta
\in
\Theta}\inf_{x \in\Dset} \g[\theta](x,Y_0)} < \infty$.
\end{longlist}
This condition is satisfied if $(x,\theta) \mapsto\g[\theta](x,y)$ is
continuous and $\Dset$ is a compact small set for all $\theta\in \Theta
$, there exists a probability measure $\nu^\theta$ such that
$\nu^\theta (\Dset)=1$ and a constant $\delta> 0$, such that, for all
$x \in \Dset $ and $A \in\Xsigma$, $\Q[\theta](x,A)
\geq\delta\nu^\theta(A)$. Note, however, that
\ref{assumlikelihoodDrift}\ref{itemcondition-minoration} is far
weaker than imposing that the set $\Dset$ is $1$-small. This is
important to deal with examples for which the transition kernel
$\Q[\theta](x,\cdot)$ does not admit a density with respect to to some
fixed dominating measure; see, for example, Section \ref
{secpostcompkalman}.
\end{rem}
%
\begin{rem}
Assumption~\ref{assumcontinuity} is\vspace*{1pt} in general the consequence of the continuity
of the kernel $\theta\mapsto\Q[\theta](x,\cdot)$ and
of the function $\theta\mapsto\g[\theta](x,\cdot)$, using classical
techniques to deal with integrals depending on a parameter.
\end{rem}
%
\begin{rem}
According to (\ref{eqP-theta-vraisemblance}), bound (\ref
{eqbound-eta-G}) may also be rewritten in terms of the kernel $\Pblock
[\theta]{\chunk{y}{0}{r-1}}$ as
\[
\sup_{x_0 \in\Cset^c} \mathbf{L}^{\theta}\bigl\langle \chunk{y} {0} {r-1} \bigr\rangle(x_0,
\Xset) \leq\eta\sup_{x_0 \in\Xset} \mathbf{L}^{\theta}\bigl\langle \chunk{y} {0}
{r-1} \bigr\rangle(x_0,\Xset)<\infty.
\]
\end{rem}
The convergence of the relative entropy is achieved for initial
distributions belonging to a particular class of initial probability
distributions. For the integer $r$ and the set $\Dset\in\Xsigma$
defined in
\ref{assumlikelihoodDrift}, let $\mathcal{M}(\Dset,r)$ be the subset
$\mathcal{P}(\Xset,\Xsigma)$ of probability measures on $(\Xset,\Xsigma
)$ satisfying
%
\begin{eqnarray}
\label{eqmort-de-rire} \mathcal{M}(\Dset,r) &=& \Bigl\{ \chi\in\mathcal
{P}(\Xset,
\Xsigma),
\nonumber\\[-8pt]\\[-8pt]
&&\hspace*{6.5pt}\mathbb E\Bigl[{\ln^- \inf_{\theta\in\Theta}
\chi\mathbf{L}^{\theta}\bigl\langle \chunk{Y} {0} {u-1} \bigr\rangle
\one_{\Dset}}\Bigr] < \infty \mbox{ for all $u \in\{1,\ldots,r\} $}
\Bigr\}. \nonumber
\end{eqnarray}

\begin{prop} \label{proplimiteRelativeEntropyRateY}
Assume~\ref{assumlikelihoodDrift} and~\ref{assummajo-g}. Then:

\renewcommand\thelonglist{(\roman{longlist})}
\renewcommand\labellonglist{\thelonglist}
\begin{longlist}
\item\label{itemlim-ponct-cond-rate-entropy-Y} for any $\theta\in
\Theta$, there exists a measurable function $\pi^\theta_Y\dvtx  \Yset^{\zset
^{-}} \to\rset$ such that for any probability measure $\chi\in
\mathcal{M}(\Dset,r)$,
\[
\mathbb P\Bigl[{\lim_{m \to\infty} p_{\chi}^{\theta}\bigl(Y_0|\chunk{Y} {-m}
{-1}\bigr)=\pi^\theta_Y\bigl(\chunk{Y} {-\infty} {0}\bigr)}\Bigr]=1;
\]
moreover,
%
\begin{equation}
\label{eqesperance-log-stat} \mathbb E\bigl[{\bigl|\ln\pi^\theta_Y\bigl(\chunk{Y}
{-\infty} {0}\bigr)\bigr|}\bigr]<\infty;
\end{equation}
\item\label{itemlim-ponct-rate-entropy-Y} for any $\theta\in\Theta$
and any probability measure $\chi\in\mathcal{M}(\Dset,r)$,
\[
\lim_{n \to\infty} n^{-1} \ln p_{\chi}^{\theta}\bigl(Y_{0}^{n-1}\bigr)
= \ell(\theta), \qquad\PP{}\mbox{-a.s.},
\]
where $\ell(\theta) \triangleq\PE{\ln\pi^\theta_Y(\chunk
{Y}{-\infty}{0})}$.
\end{longlist}
\end{prop}
%
%
\begin{theorem} \label{thmconsistanceY}
Assume~\ref{assumlikelihoodDrift}--\ref{assumcontinuity}. Then,
$\theta\mapsto\ell(\theta)$ is upper semi-continu\-ous and defining
$\Theta^\star\subset\Theta$ by $\Theta^\star\triangleq\argmax_{\theta
\in\Theta} \ell(\theta)$, we have for any probability measure $\chi
\in\mathcal{M}(\Dset,r)$,
\[
\lim_{n \to\infty} d\bigl(\mle{\chi} {n},\Theta^\star\bigr)=0,\qquad \PP{}
\mbox{-a.s.}
\]
\end{theorem}
%
\begin{rem}
When the model is well specified, the law of the observations belongs
to the parametric family of distributions on which the maximization
occurs and is therefore associated to a specific parameter $\theta^*$.
In this particular case, under some appropriate assumptions, the set
$\Theta^*$ is reduced to the singleton $\{\theta^*\}$, and the
consistency result of the MLE in well specified models can then be
written as (see~\cite{doucmoulinesolssonvanhandel2011})
\[
\lim_{n \to\infty} d\bigl(\mle{\chi} {n},\theta^\star\bigr)=0,\qquad \PP{}
\mbox{-a.s.}
\]
\end{rem}
A simple sufficient condition can be proposed to ensure that $\chi
\in\mathcal{M}(\Dset,r)$.

\begin{prop}
\label{propsufficient-condition-mdr}
Assume there exist a sequence of sets $\Dset_u \in\Xsigma$, $u \in\{
0,\ldots, r-1\}$, such that (setting $\Dset_r= \Dset$ for notational
convenience), for some $\delta> 0$,
%
\begin{equation}
\label{eqCS2-mort-de-rire} \inf_{x_{u-1} \in\Dset_{u-1}} \inf_{\theta\in
\Theta} \Q[\theta
](x_{u-1},\Dset_u) \geq\delta,\qquad u \in\{1, \ldots, r\},
\end{equation}
and
%
\begin{equation}
\label{eqCS1-mort-de-rire} \mathbb E\Bigl[{\ln^- \inf_{\theta\in\Theta} \inf_{x \in
\Dset_u}
\g[\theta](x,Y_0)}\Bigr] < \infty\qquad \mbox{for $u
\in\{0,\ldots,r\}$}.
\end{equation}
Then, any initial distribution $\chi$ on $(\Xset,\Xsigma)$ satisfying
$\chi(\Dset_0) > 0$ belongs to $\mathcal{M}(\Dset,r)$.
\end{prop}
%
\begin{rem}
\label{remcheck-CS1-mort-de-rire}
To check (\ref{eqCS1-mort-de-rire}), we typically assume that, for any
given $y \in\Yset$,
the function $(x,\theta) \mapsto\g[\theta](x,y)$ is continuous and
that $\Dset_i \times\Theta$ is a compact set, $i \in\{0,\ldots,r-1\}$.
This condition then translates into an assumption on some generalized
moments of the process $Y$.
\end{rem}

To check (\ref{eqCS2-mort-de-rire}), the following lemma is useful.
%
\begin{lem}
\label{lemcheck-CS2-mort-de-rire} Assume that $\Xset= \rset^d$ for
some integer $d > 0$ and that $\Xsigma$ is the associated Borel
$\sigma$-field. Assume in addition that, for any open subset
$\Oset\in\Xsigma$, the function $(x,\theta) \to\Q[\theta](x,\Oset)$ is
lower semi-continuous on the product space $\Xset\times\Theta$. Then,
for any $\delta> 0$ and any compact subset $\Dset_0 \in\Xsigma$, there
exists a sequence of compact subsets $\Dset_u$, $u \in\{0,\ldots,r-1\}$
satisfying (\ref{eqCS2-mort-de-rire}).
\end{lem}
%
\section{Applications}
\label{secapplications}
In this section, we develop three classes of examples. In Section
\ref{secpostcompkalman} we consider linear Gaussian state space
models. This is obviously a very important model,
which is used routinely to analyze time-series models. We analyze this
model under assumptions which are\vadjust{\goodbreak} very general and might serve to
illustrate the stated assumptions.
In Section~\ref{secfinite-state-leroux}, we consider the classic case
where state space of the underlying Markov chain is a finite set.
Finally, in Section~\ref{secexample-nonlinear-state-space}, we develop
a general class
of nonlinear state space models. In all these examples, we will find
that the assumptions of Theorem~\ref{thmconsistanceY} are satisfied
under general assumptions.

\subsection{Gaussian linear state space models}
\label{secpostcompkalman}
Gaussian linear state space models form an important class of HMMs. In
this setting, let $\Xset=\rset^{d_x}$, and $\Yset=\rset^{d_y}$ for
some integers
and let $\Theta$ be a compact parameter space. The model is specified
by
%
\begin{eqnarray}
\label{eqlinearstatespacetimedependdynamic}
X_{k+1} & = & A_{\theta} X_k + \URoot[\theta]
U_{k},
\\
\label{eqlinearstatespacetimedependobservation}
Y_k & = & B_{\theta} X_k + \VRoot[\theta]
V_k,
\end{eqnarray}
where $\{(U_k,V_k)\}_{k \geq0}$ is an i.i.d. sequence of Gaussian
vectors with zero mean and identity covariance matrix, independent of $X_0$.
Here $U_k$ is $d_u$-dimensional, $V_k$ is $d_y$-dimensional and the matrices
$A_\theta,R_\theta,B_\theta,S_\theta$ have the appropriate dimensions.

For any integer $n$, define $\mathcal{O}_{\theta,n}$ and $\mathcal
{C}_{\theta,n}$
the observability matrix and the controllability matrices
%
\begin{equation}
\label{eqdefinition-controlabilite-observabilite} \mathcal{O}_{\theta,n}
\triangleq\lleft[\matrix{ B_{\theta}
\cr
B_{\theta}
A_{\theta}
\cr
B_{\theta} A_{\theta}^2
\cr
\vdots
\cr
B_{\theta} A_{\theta}^{n-1} } \rright] \quad\mbox{and}\quad
\mathcal{C}_{\theta,n}\triangleq\bigl[\matrix{ A_\theta^{n-1}
\URoot[\theta] A_\theta^{n-2} \URoot[\theta] \cdots\URoot[
\theta] } \bigr].
\end{equation}
It is assumed in the sequel that for any $\theta\in\Theta$, the following
hold:
\renewcommand\thelonglist{(L\arabic{longlist})}
\renewcommand\labellonglist{\thelonglist}
\begin{longlist}
\item\label{hypLGSScommandable-observable}
The pair $[A_\theta,B_\theta]$ is observable, and the
pair $[A_\theta,\URoot[\theta]]$ is controllable; that is, there exists
an integer $r$ such that, the observability
matrix $\mathcal{O}_{\theta,r}$ and the controllability matrix
$\mathcal{C}_{\theta,r}$ are full rank.
\end{longlist}
\begin{longlist}
\setcounter{longlist}{1}
\item
\label{hypLGSSobservation-noise-full-rank}
The measurement noise covariance matrix $\VRoot[\theta]$ is full rank.
\end{longlist}
\begin{longlist}
\setcounter{longlist}{2}
\item
\label{hypLGSScontinuity-matrix}
The functions $\theta\mapsto A_{\theta}$, $\theta\mapsto\URoot
[\theta
]$, $\theta\mapsto B_{\theta}$ and $\theta\mapsto
\VRoot[\theta]$ are continuous on~$\Theta$.
\end{longlist}
\begin{longlist}
\setcounter{longlist}{3}
\item\label{hypLGSSvariance-observation} $\PE{\|Y_0\|^2} < \infty$.
\end{longlist}

We now check the assumptions of Theorem~\ref{thmconsistanceY}.

The dimension $d_u$ of the state noise vector $U_k$ is in many situations
smaller than the dimension $d_x$ of the state vector $X_k$ and hence
$\UCov[\theta]$ (where ${}^tA$ is the transpose of the matrix $A$)
may be rank deficient.

Some additional notation is needed.
For any positive matrix $A$ and any vector $z$ of appropriate
dimension, denote $\| z \|_A^2= {}^tz A^{-1} z$.
Define for any integer $n$
%
\begin{equation}
\label{eqdefinition-f} \mathcal{F}_{\theta,n} = \mathcal{D}_{\theta,n}
{}^t\mathcal{D}_{\theta,n} + \mathcal{S}_{\theta,n}
{}^t\mathcal{S}_{\theta,n},
\end{equation}
where ${}^t$ denotes the transpose and
\begin{eqnarray*}
\mathcal{D}_{\theta,n} &\triangleq& \lleft[ \matrix{ 0 & 0 & \cdots& 0
\cr
B_\theta R_\theta& \ddots& & 0
\cr
B_\theta
A_\theta R_\theta& B_\theta R_\theta& \ddots&
\vdots
\cr
\vdots& & \ddots& 0
\cr
B_\theta A_\theta^{n-2}
R_\theta& B_\theta A_\theta^{n-3}
R_\theta& \cdots& B_\theta R_\theta} \rright],
\\
\mathcal{S}_{\theta,n}&\triangleq&\lleft[ \matrix{ S_\theta& 0 &
\cdots& 0
\cr
0 & S_\theta& \ddots& \vdots
\cr
\vdots& \ddots& \ddots& 0
\cr
0 & \cdots& 0& S_\theta} \rright].
\end{eqnarray*}
Under~\ref{hypLGSSobservation-noise-full-rank}, for any $n \geq
r$, the matrix $\mathcal{F}_{\theta,n}$ is positive definite.
The likelihood of the observations $\chunk{y}{0}{n-1} \in\Yset^n$
starting from $x_0$ is given by
%
\begin{equation}
\label{eqLGSSlikelihood-chunk-observation}
p_{x_0}^{\theta}\bigl(y_{0}^{n-1}\bigr)
= (2 \pi)^{-n d_y} \operatorname{det}^{-1/2}(
\mathcal{F}_{\theta,n}) \exp\bigl( - \tfrac12 \llVert\mathbf{y}_{n-1}
- \mathcal{O}_{\theta,n} x_0 \rrVert^2_{\mathcal
{F}_{\theta,n}}
\bigr),
\end{equation}
where $\mathbf{y}_{n-1}= {}^t{[{}^ty_0, {}^ty_1,\ldots,
{}^ty_{n-1}]}$, and $\mathcal{O}_{\theta,n}$ is defined in (\ref
{eqdefinition-controlabilite-observabilite}).

Consider first~\ref{assumlikelihoodDrift}.
Under~\ref{hypLGSScommandable-observable}, the observability
matrix $\mathcal{O}_{\theta,r}$ is full rank,
we have, for any compact subset $\Kset\subset\Yset^{r}$,
\[
\lim_{\| x_0 \| \to\infty} \inf_{\chunk{y}{0}{r-1} \in\Kset} \llVert
\mathbf{y}_{r-1} -
\mathcal{O}_{\theta,r} x_0 \rrVert_{\mathcal
{F}_{\theta,r}} = \infty,
\]
showing that, for all $\eta> 0$, we may choose a compact set $\Cset$
in such a way that (\ref{eqbound-eta-G}) is satisfied. It remains to
prove that any compact set $\Cset$ is a $r$-local Doeblin satisfying
the condition (\ref{eqlower-bound}).
For any $\chunk{y}{0}{r-1} \in\Yset^{r-1}$ and $x_0 \in\Xset$
the measure $\Pblock[\theta]{\chunk{y}{0}{r-1}}(x_0,\cdot)$ is
absolutely continuous
with respect to the Lebesgue measure on $\Xset$ with Radon--Nikodym
denoted $\pblock[\theta]{\chunk{y}{0}{r-1}}(x_0,x_r)$
given (up to an irrelevant multiplicative factor) by
%
\begin{equation}
\label{eqGLSSradon-nikodym}\quad \bolds{\ell}^{\theta}\bigl\langle\chunk{y} {0}
{r-1} \bigr\rangle (x_0,x_r) \propto\operatorname{det}^{-1/2}(
\mathcal{G}_{\theta,r}) \exp\biggl( - \tfrac {1}{2} \biggl\llVert\lleft[
\matrix{ \mathbf{y}_{r-1} \cr x_r } \rright] - \lleft[ \matrix{
\mathcal{O}_{\theta,r} \cr A_\theta^r } \rright] x_0\biggr
\rrVert_{\mathcal {G}_{\theta,r}}^2 \biggr),\hspace*{-22pt}
\end{equation}
where the covariance matrix $\mathcal{G}_{\theta,r}$ is given by
\[
\mathcal{G}_{\theta,r} = \lleft[\matrix{ \mathcal{D}_{\theta,r}
\cr
\mathcal{C}_{\theta,r}} \rright] \bigl[ \matrix{ {}^t
\mathcal{D}_{\theta,r} {}^t\mathcal{C}_{\theta,r}} \bigr] + \lleft[
\matrix{ \mathcal{S}_{\theta,r}
\cr
\mathbf{0} } \rright]
\bigl[ \matrix{ {}^t\mathcal{S}_{\theta,r} {}^t
\mathbf{0} } \bigr].
\]
The proof of (\ref{eqGLSSradon-nikodym}) relies on the positivity of
$\mathcal{G}_{\theta,r}$, which requires further discussion.
By construction, the matrix $\mathcal{G}_{\theta,r}$ is nonnegative.
For any $\mathbf{y}_{r-1} \in\Yset^r$ and $x \in\Xset$, the equation
\[
\bigl[{}^t\mathbf{y}_{r-1} {}^tx\bigr]
\mathcal{G}_{\theta,r} \lleft[ \matrix{ \mathbf{y}_{r-1}
\cr
x }
\rright] = \bigl\llVert{}^t\mathcal{D}_{\theta,r}
\mathbf{y}_{r-1}+ {}^t\mathcal{C}_{\theta,r} x \bigr
\rrVert^2 + \bigl\llVert{}^t\mathcal{S}_{\theta,r}
\mathbf{y}_{r-1} \bigr\rrVert^2 = 0
\]
implies that $\llVert {}^t\mathcal{D}_{\theta,r} \mathbf{y}_{r-1}+
{}^t\mathcal{C}_{\theta,r} x \rrVert^2 = 0$ and
$\llVert {}^t\mathcal{S}_{\theta,r} \mathbf{y}_{r-1} \rrVert^2 =0
$. Since the matrix $\mathcal{S}_{\theta,r}$ is full rank, this
implies that
$\mathbf{y}_{r-1}= 0$. Since $\mathcal{C}_{\theta,r}$ is full-rank (the
pair $[A_\theta,R_\theta]$ is controllable),
this implies that $x=0$. Therefore, the matrix $\mathcal{G}_{\theta,r}$
is positive definite and, for any $\mathbf{y}_{r-1}$,
the function
\[
(x_0,x_r) \mapsto\biggl\llVert\lleft[ \matrix{
\mathbf{y}_{r-1}
\cr
x_r } \rright] - \lleft[ \matrix{
\mathcal{O}_{\theta,r}
\cr
A_\theta^r } \rright]
x_0\biggr\rrVert_{\mathcal{G}_{\theta,r}}^2
\]
is continuous, and is therefore bounded on any compact subset of $\Xset
\times\Xset$.
This implies that every nonempty compact set $\Cset\subset\rset^{d_x}$
is a $r$-local Doeblin set, with $\lambda_\Cset^\theta(\cdot
) =
\lambda^{\mathrm{Leb}}(\cdot)/ \lambda^{\mathrm{Leb}}(\Cset)$ and
\begin{eqnarray*}
\epsilon^-_\Cset\bigl(\chunk{y} {0} {r-1}\bigr) &=& \bigl( \lambda^{\mathrm
{Leb}}(
\Cset) \bigr)^{-1} \inf_{\theta\in\Theta} \inf_{(x_0,x_r) \in\Cset\times
\Cset} \bolds{\ell}^{\theta}\bigl\langle \chunk{y} {0} {r-1} \bigr\rangle
(x_0,x_r),
\\
\epsilon^+_\Cset\bigl(\chunk{y} {0} {r-1}\bigr) &=& \bigl( \lambda^{\mathrm
{Leb}}(
\Cset) \bigr)^{-1} \sup_{\theta\in\Theta} \sup_{(x_0,x_r) \in\Cset\times
\Cset} \bolds{\ell}^{\theta}\bigl\langle \chunk{y} {0} {r-1} \bigr\rangle (x_0,x_r).
\end{eqnarray*}
Therefore, condition (\ref{eqlower-bound}) is satisfied with any
compact set $\Kset\subseteq\Yset^{r-1}$.
It remains to show~\ref{assumlikelihoodDrift}\ref
{itemcondition-minoration}.
Under~\ref{hypLGSScommandable-observable}, $\Pblock[\theta]{\chunk
{y}{0}{r-1}}(x_0,\cdot)$ is absolutely continuous with respect to the Lebesgue
measure $\lambda^{\mathrm{Leb}}$. Therefore, for any set $\Dset$,
\[
\inf_{\theta\in\Theta} \inf_{x_0 \in\Dset} \mathbf{L}^{\theta}\bigl\langle \chunk{y} {0}
{r-1} \bigr\rangle(x_0,\Dset) \geq\inf_{\theta\in\Theta} \inf_{(x_0,x_r) \in\Dset
\times\Dset}
\bolds{\ell}^{\theta}\bigl\langle \chunk{y} {0} {r-1} \bigr\rangle (x_0,x_r)
\lambda^{\mathrm{Leb}}(\Dset).
\]
Take $\Dset$ to be any compact set with
positive Lebesgue measure.
\begin{eqnarray*}
&& \sup_{\theta\in\Theta} \sup_{(x_0,x_r) \in\Dset\times\Dset} \biggl
\llVert\lleft[ \matrix{
\mathbf{y}_{r-1}
\cr
x_r } \rright] - \lleft[ \matrix{
\mathcal{O}_{\theta,r}
\cr
A_\theta^r } \rright]
x_0\biggr\rrVert_{\mathcal{G}_{\theta,r}}^2
\\
&&\qquad \leq2 \lambda_{\max} (\mathcal{G}_{\theta,r} ) \Bigl\{ \|
\mathbf{y}_{r-1} \|^2 + \max_{x \in\Dset} \|x
\|^2 \bigl[1+\lambda_{\max} \bigl({}^t
\mathcal{O}_{\theta,r} \mathcal{O}_{\theta,r} + {}^tA_\theta^r
A^r_\theta\bigr) \bigr] \Bigr\},
\end{eqnarray*}
where $\lambda_{\max}(A)$ is the largest eigenvalue of $A$. Under
\ref{hypLGSScontinuity-matrix}, $\theta\mapsto
\lambda_{\max} (\mathcal{G}_{\theta,r} )$ and $\theta
\mapsto
\lambda_{\max} ({}^t\mathcal{O}_{\theta,r}
\mathcal{O}_{\theta,r} + {}^tA_{\theta}^r A^r_{\theta}
)$ are
bounded. Under~\ref{hypLGSSvariance-observation}, $\PE{\|Y_0\|^2}
< \infty$, then
\ref{assumlikelihoodDrift}\ref{itemcondition-minoration} is
satisfied for any compact set.

Consider now~\ref{assummajo-g}. Under
\ref{hypLGSSobservation-noise-full-rank}, $\VRoot[\theta]$ is full
rank, and choosing the reference measure $\mu$ to be the Lebesgue
measure on $\Yset$, we find that $\g[\theta](x,y)$ is a Gaussian
density for each $x\in\Xset$ with covariance matrix $\VCov[\theta
]$. We
therefore have
\[
\sup_{\theta\in\Theta} \sup_{x \in\Xset} \g[\theta](x,y) = (2\pi
)^{-d_y/2} \sup_{\theta\in\Theta} \operatorname{det}^{-1/2}\bigl(\VCov[
\theta]\bigr)<\infty,
\]
so that~\ref{assummajo-g}\ref{itemassum-pos-g} and \ref
{itemassum-supg} are satisfied.

We finally check~\ref{assumcontinuity}. For any $n \geq r$, and $x
\in
\Xset$
the function $ \theta\mapsto\dens[x_0][\theta]{\chunk{y}{0}{n-1}}$ is
given by (\ref{eqLGSSlikelihood-chunk-observation}).
Under~\ref{hypLGSScontinuity-matrix}, the functions $\theta
\mapsto\mathcal{O}_{\theta,n}$ [where $\mathcal{O}_{\theta,n}$ is the
observability matrix defined in (\ref
{eqdefinition-controlabilite-observabilite})] and $\theta\mapsto
\operatorname{det}^{-1/2}(\mathcal{F}_{\theta,n})$ [where $\mathcal
{F}_{\theta,n}$ is the covariance matrix defined in (\ref
{eqdefinition-f})], are
continuous on $\Theta$ for any $n\ge r$.
Thus, for any $x \in\Xset$, $\theta\mapsto\dens[x][\theta]{\chunk
{y}{0}{n-1}}$
is continuous for every $n\ge r$, showing~\ref{assumcontinuity}.\vadjust{\goodbreak}

To conclude this discussion, we need to specify more explicitly the set
$\mathcal{M}(\Dset,r)$ [see (\ref{eqmort-de-rire})]
of possible initial distributions. Using Proposition
\ref{propsufficient-condition-mdr},
we have to check the sufficient conditions (\ref{eqCS2-mort-de-rire})
and (\ref{eqCS1-mort-de-rire}).
To check (\ref{eqCS2-mort-de-rire}), we use Lemma
\ref{lemcheck-CS2-mort-de-rire}. Note that, for any open subset $\Oset$,
\[
\Q[\theta](x,\Oset)= \mathbb E\bigl[{\one_{\Oset}(A_\theta x +
R_\theta U)}\bigr],
\]
where the expectation is taken with respect to the standard normal random
variable~$U$. Let $\{ (x_n,\theta_n) \}_{n=1}^\infty$ be any sequence
converging to $(x,\theta)$. By the Fatou lemma, using that function
$\one_{\Oset}$ is lower semi-continuous and that $\theta\mapsto A_\theta$
is continuous under~\ref{hypLGSScontinuity-matrix}, we have
\begin{eqnarray*}
\liminf_{n \to\infty} \Q[\theta_n](x_n,\Oset) &\geq&\mathbb E
\Bigl[{\liminf_{n \to
\infty} \one_{\Oset}(A_{\theta_n} x +
R_{\theta_n} U)}\Bigr]
\\
&\geq&\mathbb E
\Bigl[{\liminf_{n \to\infty} \one_{\Oset}(A_{\theta_n} x +
R_{\theta
_n} U)}\Bigr] = \Q[\theta](x,\Oset),
\end{eqnarray*}
showing that, for any open subset $\Oset$,
the function $(x,\theta) \mapsto\Q[\theta](x,\Oset)$ is lower
semi-continuous.

Assumption~\ref{hypLGSSobservation-noise-full-rank} implies that,
for all $(x,y) \in\Xset\times\Yset$,
\begin{eqnarray*}
\ln\g[\theta](x,y) &\geq&-\frac{d_y}{2} \ln( 2 \pi) -\frac12
\inf_{\theta\in\Theta} \ln\operatorname{det}^{-1/2}\bigl(\VCov[\theta]\bigr)
\\
&&{}- \Bigl[\inf_{\theta\in\Theta}\lambda_{\min} \bigl(\VCov[\theta] \bigr)
\Bigr]^{-1} \Bigl[ \|y\|^2 + \sup_{\theta\in\Theta}\|
B_\theta x\|^2 \Bigr],
\end{eqnarray*}
where $\lambda_{\min}(\VCov[\theta])$ is the minimal eigenvalue of
$\VCov[\theta]$. Therefore, under~\ref{hypLGSSvariance-observation},
(\ref{eqCS1-mort-de-rire})~is satisfied because $\Dset_u$ is a compact
set, $u \in\{0,\ldots,r\}$.

We can therefore apply Theorem~\ref{thmconsistanceY} to show that
the MLE is consistent for any initial measure $\chi$ as soon as the
process $\{ Y_k \}_{k \in\zset}$
is stationary ergodic and $\PE{|Y_0|^2} < \infty$.

\subsection{Finite state models}
\label{secfinite-state-leroux}

One of the most widely used classes of HMMs is obtained when the
state-space is finite, that is,
$\Xset=\{1,\ldots,d\}$ for some integer $d$, $\Yset$ is any Polish
space and $\Theta$ is a compact metric space. For each parameter
$\theta\in\Theta$, the transition kernel $\Q[\theta]$ is determined
by the corresponding transition probability matrix $\mathcal{Q}_\theta$,
while the observation density $\g[\theta]$ is given as in the general
setting of this paper.

It is assumed in the sequel that:
\renewcommand\thelonglist{(F\arabic{longlist})}
\renewcommand\labellonglist{\thelonglist}
\begin{longlist}
\item\label{itemFuniform-irreducibility}
There exists an integer $r\! >\! 0$, such that,
\mbox{$\inf_{\theta\in\Theta} \inf_{(x,x') \in\Xset\times\Xset}
\mathcal
{Q}^r_\theta(x,x') \!>\! 0$}.
\end{longlist}
\begin{longlist}
\setcounter{longlist}{1}
\item\label{itemFupper-and-lower-bound-likelihood}
There exists a set $\Mset\subset\Yset$ such that
$
\inf_{\theta\in\Theta} \inf_{y \in\Mset} \inf_{x \in\Xset} \g
[\theta
](x,y) > 0
$
and
$
\sup_{\theta\in\Theta} \sup_{y \in\Mset} \sup_{x \in\Xset} \g
[\theta
](x,y) < \infty.
$
\end{longlist}
\begin{longlist}
\setcounter{longlist}{2}
\item\label{itemFmajo-g}
For any $\theta\in\Theta$, the function $\g[\theta]\dvtx  (x,y) \in
\Xset
\times\Yset\mapsto\g[\theta](x,y)$ is positive and
\[
\mathbb E
\Bigl[{ \ln^+\sup_{\theta\in\Theta} \sup_{x \in\Xset} \g[\theta
](x,Y_0)}\Bigr]<\infty.\vadjust{\goodbreak}
\]
\end{longlist}
\begin{longlist}
\setcounter{longlist}{3}
\item\label{itemFmino-g}
$\PE{\ln^- \inf_{\theta\in\Theta} \inf_{x \in\Xset} \g[\theta
](x,Y_0)} < \infty$.
\end{longlist}
\begin{longlist}
\setcounter{longlist}{4}
\item\label{itemFcontinuity}
$\theta\mapsto\mathcal{Q}_\theta$ and $\theta\mapsto g_\theta(x,y)$
are continuous for any $x\in\Xset$, $y\in\Yset$.
\end{longlist}

Consider first~\ref{assumlikelihoodDrift}. We set $\Cset= \Xset$.
Since $\Cset^c=\varnothing$, (\ref{eqbound-eta-G}) is trivially
satisfied. Under~\ref{itemFuniform-irreducibility},
equation (\ref{eqdefinition-LD-set}) is satisfied with
$\varphi_{\Xset} \langle\chunk{y}{0}{r-1}\rangle(x) \equiv1$,
$\lambda_\Xset^\theta= d^{-1} \sum_{i=1}^{d} \delta_{i}$, and
\begin{eqnarray*}
\epsilon_\Xset^{-} \bigl[\chunk{y} {0} {r-1}\bigr] &=& d \prod
_{i=0}^{d-1} \inf_{\theta\in\Theta}
\inf_{x \in\Xset} \g[\theta](x,y_i) \times\inf_{\theta\in\Theta}
\inf_{(x,x') \in\Xset\times\Xset} \mathcal{Q}_\theta^r
\bigl(x,x'\bigr),
\\
\epsilon_\Xset^+ \bigl[\chunk{y} {0} {r-1}\bigr] &=& d \prod
_{i=0}^{d-1} \sup_{\theta\in\Theta} \sup_{x \in\Xset}
\g[\theta](x,y_i) \times\sup_{\theta\in\Theta} \sup_{(x,x') \in\Xset
\times\Xset}
\mathcal{Q}_\theta^r\bigl(x,x'\bigr).
\end{eqnarray*}
Hence, the\vspace*{1pt} state space $\Xset$ is a $r$-local Doeblin set. Assumption
\ref{itemFupper-and-lower-bound-likelihood} implies that (\ref
{eqlower-bound}) is satisfied with $\Kset= \Mset^r$. Now, note that for
all $u \in\{1, \ldots, r\}$ and $\chunk{y}{0}{u-1} \in\Yset^r$,
%
\begin{equation}
\label{eqlower-bound-Pblock} \inf_{\theta\in\Theta} \inf_{x \in\Xset}
\mathbf{L}^{\theta}\bigl\langle \chunk{y} {0} {u-1} \bigr\rangle \geq\prod_{i=0}^{u-1}
\inf_{\theta\in\Theta} \inf_{x
\in\Xset} \g[\theta](x,y_i).
\end{equation}
Using the previous inequality with $u=r$ and noting that
\ref{itemFmino-g} implies that $\PE{\ln^- \inf_{\theta\in\Theta
} \inf_{x \in\Xset} \g[\theta](x,Y_0)} < \infty$, we see that
equation (\ref{eqcondition-minoration}) is satisfied with $\Dset=
\Xset
$. The same argument for any $u\in{1,\ldots,r}$ shows that all the
probability measures on $(\Xset,\Xsigma)$ belong to the set $\mathcal
{M}(\Xset,r)$, defined in (\ref{eqmort-de-rire}).

Assumption~\ref{assummajo-g} is a direct consequence of
\ref{itemFmajo-g}.
Finally, we note that the continuity of $\theta\mapsto\mathcal
{Q}_\theta
$ and $\theta\mapsto g_\theta(x,y)$
yield immediately that $\theta\mapsto\dens[x][\theta]{\chunk
{y}{0}{n}}$ is
a continuous function for every $n\ge0$ and $\chunk{y}{0}{n} \in
\Yset^{n+1}$,
establishing~\ref{assumcontinuity}.

We can therefore apply Theorem~\ref{thmconsistanceY} under
\ref{itemFuniform-irreducibility}--\ref{itemFcontinuity} to show that
the MLE is consistent for any initial measure $\chi$ as soon as the
process $\{ Y_k \}_{k \in\zset}$
is stationary ergodic.

\subsection{Nonlinear state space models}
\label{secexample-nonlinear-state-space}

In this section, we consider a class of nonlinear state space models.
Let $\Xset=\rset^d$, $\Yset=\rset^\ell$ and $\Xsigma$ and
$\Ysigma$ be the
associated Borel $\sigma$-fields. Let $\Theta$ be a compact
metric space. For each $\theta\in\Theta$ and each $x \in\Xset$,
the Markov kernel $Q_\theta(x,\cdot)$ has a density $q_\theta
(x,\cdot)$
with respect to the Lebesgue measure on $\Xset$.

For example, $(X_k)_{k \ge0}$ may be defined through the nonlinear recursion
\[
X_{k} = T_\theta(X_{k-1})+ \Sigma_\theta(X_{k-1})
\zeta_k,
\]
where $(\zeta_k)_{k \geq1}$ is an i.i.d. sequence of $d$-dimensional
random vectors which are assumed to possess a density $\rho_\zeta$ with
respect to the Lebesgue measure $\lambda^{\mathrm{Leb}}$ on $\rset^d$,
and $T_\theta\dvtx
\rset^d \to\rset^d$, $\Sigma_\theta\dvtx  \rset^d \to\rset^{d \times d}$
are given (measurable) matrix-valued functions such that for each
$\theta\in\Theta$ and $x \in\Xset$,
$\Sigma_\theta(x)$ is full-rank. Such a model for
$(X_k)_{k \geq0}$ is sometimes known as a vector ARCH model, and covers
many models\vadjust{\goodbreak} of interest in time series analysis and financial
econometrics. We let the reference measure $\mu$ be the Lebesgue
measure on
$\rset^\ell$, and define the observed process $(Y_k)_{k\ge0}$ by means
of a given observation density $\g[\theta](x,y)$.

We now introduce the basic assumptions of this section.

\renewcommand\thelonglist{(NL\arabic{longlist})}
\renewcommand\labellonglist{\thelonglist}
\begin{longlist}
\item\label{itemNLpositivity-density}
The function $ (x,x',\theta) \mapsto\q[\theta](x,x')$ is a positive
continuous function on $\Xset\times\Xset\times\Theta$. In addition,
$\sup_{\theta\in\Theta} \sup_{(x,x') \in\Xset\times\Xset} \q
[\theta
](x,x') < \infty$.
\end{longlist}
\begin{longlist}
\setcounter{longlist}{1}
\item\label{itemNLgoes-to-zero-at-infinity}
For any compact subset $\Kset\subset\Yset$, and $\theta\in\Theta$,
\[
\lim_{|x|\to\infty} \sup_{y \in\Kset} \frac{\g[\theta
](x,y)}{\sup_{x'
\in\Xset} \g[\theta](x',y)} = 0.
\]
\end{longlist}
\begin{longlist}
\setcounter{longlist}{2}
\item\label{itemNLbound-log+-g}
For each $(x,y) \in\Xset\to\Yset$, the function $\theta\mapsto\g
[\theta](x,y)$ is positive and continuous on $\Theta$. Moreover,
\[
\mathbb E
\Bigl[{\ln^+ \sup_{\theta\in\Theta} \sup_{x\in\Xset} \g[\theta
](x,Y_0) }\Bigr] <\infty.
\]
\end{longlist}
\begin{longlist}
\setcounter{longlist}{3}
\item\label{itemNLbound-log--g}
There exists a compact subset $\Dset\subset\Yset$ such that
\[
\mathbb E
\Bigl[{ \ln^- \inf_{\theta\in\Theta} \inf_{x \in\Dset} \g[\theta
](x,Y_0)}\Bigr] < \infty.
\]
\end{longlist}
We have made no attempt at generality here: for sake of example, we have
chosen a set of conditions under which the assumptions of Theorem
\ref{thmconsistanceY} are easily verified. Of course, the applicability
of Theorem~\ref{thmconsistanceY} extends far beyond the simple
assumptions imposed in this section.
%
\begin{rem}
Nonetheless, the present assumptions already cover a broad class of
nonlinear models. Consider, for example, the stochastic volatility
model~\cite{hullwhite1987}
%
\begin{equation}
\cases{ X_{k+1}= \phi_\theta X_k
+ \sigma_\theta\zeta_k,
\cr
Y_k= \beta_\theta\exp( X_k/2)
\varepsilon_k,}
\end{equation}
where $(\zeta_k,\varepsilon_k)$ are i.i.d. Gaussian random variables in
$\rset^2$ with zero mean and identity covariance matrix, $\beta_\theta>
0$, $\sigma_\theta> 0$ for every $\theta\in\Theta$, and the functions
$\theta\mapsto\phi_\theta$, $\theta\mapsto\sigma_\theta$, and
$\theta\mapsto\beta_\theta$ are continuous. Then, assumptions
\ref{itemNLpositivity-density}--\ref{itemNLbound-log--g} are satisfied
as noted by Douc et al.~\cite{doucmoulinesolssonvanhandel2011},
Remark 10.

%
%
\end{rem}
Under~\ref{itemNLpositivity-density}, every compact set $\Cset
\subset\Xset= \rset^d$ with
$\lambda^{\mathrm{Leb}}(\Cset) > 0$ is a 1-small set and therefore a
local Doeblin with
$\lambda^\theta_\Cset(\cdot)= \lambda^{\mathrm{Leb}}(\cdot\cap
\Cset)/ \lambda^{\mathrm{Leb}}(\Cset)$,
$\varphi_{\Cset}^\theta\langle y_0\rangle=\lambda^{\mathrm
{Leb}}(\Cset) $ and
\begin{eqnarray*}
\epsilon_\Cset^- &=& \inf_{\theta\in\Theta} \inf_{(x,x')\in\Cset
\times\Cset} \q[
\theta]\bigl(x,x'\bigr),
\\
\epsilon_\Cset^+ &=& \sup_{\theta\in\Theta} \sup_{(x,x') \in
\Cset
\times\Cset} \q[
\theta]\bigl(x,x'\bigr).
\end{eqnarray*}
Under~\ref{itemNLpositivity-density} and
\ref{itemNLgoes-to-zero-at-infinity}, (\ref{eqbound-eta-G}) and
(\ref{eqlower-bound}) are satisfied with $r=1$; equation (\ref
{eqcondition-minoration}) follows from~\ref{itemNLpositivity-density}
and~\ref{itemNLbound-log--g}. Thus assumption
\ref{assumlikelihoodDrift} holds.\vadjust{\goodbreak}

Assumption
\ref{assummajo-g} follows directly from~\ref{itemNLbound-log+-g}.
To establish~\ref{assumcontinuity}, it suffices to note that, under
\ref{itemNLpositivity-density}, for any $(x,x') \in\Xset\times
\Xset$, $\theta\mapsto\q[\theta](x,x')$ is continuous, under
\ref{itemNLbound-log+-g}, for any $(x,y) \in\Xset\times\Yset$,
$\theta\mapsto\g[\theta](x,y)$ is continuous and for any \mbox{$n \in
\nset$},
$\sup_{\theta\in\Theta} \sup_{x \in\Xset} \prod_{k=0}^n \g
[\theta
](x,Y_k) < \infty$, $\PP{}$-a.s.
The bounded convergence theorem shows that, $\PP{}$-a.s. the function
$\theta\mapsto\dens[x][\theta]{\chunk{Y}{0}{n}}$ is continuous.

Finally, under~\ref{itemNLpositivity-density}--\ref{itemNLbound-log--g}
according to Theorem~\ref{thmconsistanceY} and Proposition
\ref{propsufficient-condition-mdr} the MLE is consistent for any
initial measure $\chi$ such that $\chi(\Dset)>0$.

\section{\texorpdfstring{Proofs of Proposition \protect\ref{proplimiteRelativeEntropyRateY} and Theorem \protect\ref{thmconsistanceY}}
{Proofs of Proposition 1 and Theorem 2}}
\label{secproof}
\subsection{Block decomposition}
The first step of the proof consists of splitting the observations into
blocks of size $r$ where $r$ is defined in~\ref{assumlikelihoodDrift}.
More precisely, we will first show the equivalent of Proposition \ref
{proplimiteRelativeEntropyRateY}
and Theorem~\ref{thmconsistanceY} with $Y_i$ replaced by $Z_i
\triangleq
\chunk{Y}{ir}{(i+1)r-1}$.
With this notation,
\[
\mle{\chi} {nr}= \argmax_{\theta\in\Theta} \ln p_{\chi}^{\theta}\bigl(Y_{0}^{nr-1}\bigr)
= \argmax_{\theta\in\Theta} \ln p_{\chi}^{\theta}\bigl(Z_{0}^{n-1}\bigr).
\]
In the following, $\mle{\chi}{nr}$ is called the \textit{block maximum
likelihood estimator} (denoted hereafter as the block MLE) associated
to the observations $Z_0,\ldots,Z_{n-1}$.

\subsubsection{Forgetting of the initial distribution for the block
conditional likelihood}
Denote, for $i \in\zset$,
%
\begin{equation}
\label{eqdefinition-z} z_i=\chunk{y} {ir} {(i+1)r-1} \in
\Yset^r.
\end{equation}
Then, the likelihood $\dens[\chi][\theta]{\chunk{z}{0}{n-1}}$ may be
rewritten as
%
\begin{equation}
\label{eqvrais-en-fonction-de-P} p_{\chi}^{\theta}\bigl(z_{0}^{n-1}\bigr)=
p_{\chi}^{\theta}\bigl(y_{0}^{nr-1}\bigr)=\chi\Pblock[\theta]{z_0}
\cdots\Pblock[\theta]{z_{n-1}}\one_\Xset=\chi\mathbf{L}^{\theta}\bigl\langle
\chunk{z} {0} {n-1} \bigr\rangle\one_\Xset,
\end{equation}
where $\Pblock[\theta]{\chunk{z}{0}{n-1}}=\Pblock[\theta]{\chunk
{y}{0}{n r-1}}$ is defined in (\ref{eqdef-Pblock}).

For any sequence $ \{ z_i \}_{i \geq0} \in\Zset^{\nset}$
where $\Zset\triangleq\Yset^{r}$,
any probability measures $\chi$ and $\chi'$ on $(\Xset,\Xsigma)$
and any measurable nonnegative functions $f$ and $h$ from $\Xset$ to
$\rset^+$, define
%
\begin{eqnarray}
\label{eqdefinition-Delta} \Delta^{\theta}_{\chi,\chi'}\bigl\langle \chunk{z} {0}
{n-1} \bigr\rangle  ({f,h}) &=& \bigl(\chi\mathbf{L}^{\theta}\bigl\langle \chunk{z} {0}
{n-1} \bigr\rangle f\bigr ) \bigl( \chi'
\mathbf{L}^{\theta}\bigl\langle \chunk{z} {0} {n-1} \bigr\rangle
h\bigr)\nonumber\\[-8pt]\\[-8pt]
&&{} - \bigl(\chi \mathbf{L}^{\theta}\bigl\langle \chunk{z} {0}
{n-1} \bigr\rangle h\bigr) \bigl(\chi' \mathbf{L}^{\theta}\bigl\langle \chunk{z} {0}
{n-1} \bigr\rangle f\bigr).\nonumber
\end{eqnarray}
%
Let $\Xsetprod=\Xset\times\Xset$ and $\Xsigmaprod= \Xsigma\otimes
\Xsigma$. For $P$ a (possibly unnormalized) kernel on
$(\Xset,\Xsigma)$, we denote by $\Pprod$ the transition kernel on
$(\Xsetprod,\Xsigmaprod)$ defined, for any $(x,x') \in\Xsetprod$ and
$\Aset$, $\Aset'
\in\Xsigma$, by
%
\begin{equation}
\label{eqproduct-kernel} \Pprod\bigl[ \bigl(x,x'\bigr), \Aset
\times\Aset'\bigr]= P(x,\Aset) P\bigl(x',
\Aset'\bigr).
\end{equation}
If $\chi$ and $\chi'$ are two probability measures on $(\Xset,\Xsigma)$
and $f,g$ are real valued measurable functions on $(\Xset, \Xsigma)$,
define for $\bar{\Aset} \in\Xsigmaprod$ and $\bar{w}=(w,w')
\in\bar{\Xset}$,
%
\begin{equation}
\label{eqnotationtimes}\qquad \chi\otimes\chi'(\bar{\Aset})= \iint
\chi(\rmd x) \chi'\bigl(\rmd x'\bigr)
\one_{\bar{\Aset}}\bigl(x,x'\bigr),\qquad f\otimes h (\bar{w})=f(w)g
\bigl(w'\bigr).
\end{equation}
With the notation introduced above, (\ref{eqdefinition-Delta}) can be
rewritten as follows:
%
\begin{eqnarray}
\label{eqkey-identity} \Delta^{\theta}_{\chi,\chi'}\bigl\langle \chunk{z} {0}
{n-1} \bigr\rangle  ({f,h})&=& \int\cdots\int\chi\otimes\chi'\bigl(\rmd
\bar{w}'_0\bigr)
\Biggl(\prod_{i=0}^{n-1}
\Pblockprod[\theta]{z_i}(\bar{w}_{i},\rmd
\bar{w}_{i+1}) \Biggr)\nonumber\\[-8pt]\\[-8pt]
&&\hspace*{0pt}{}\times \{ f \otimes h - h \otimes f \}(
\bar{w}_{n}).
\nonumber
\end{eqnarray}

The following proposition extends~\cite{doucfortmoulinespriouret2009},
Proposition 12.
%
\begin{prop}
\label{propmajo-difference-numerateur}
Assume~\ref{assumlikelihoodDrift}. Let $0 \leq\gamma^- < \gamma^+
\leq1$.
Then, for any $\eta> 0$, there exists $\rho\in(0,1)$ such that,
for any sequence $(z_i)_{i \geq0} \in\Zset^{\nset}$ satisfying
%
\begin{equation}
\label{eqdef-K} n^{-1}\sum_{i=0}^{n-1}
\one_{\Kset}(z_i) \geq\max\bigl(1-\gamma^-,\bigl(1+\gamma^+
\bigr)/2 \bigr)
\end{equation}
for any $\beta\in\ooint{\gamma^-,\gamma^+}$, any nonnegative bounded
functions $f$ and $h$,
any probability measures $\chi$ and $\chi'$ on $(\Xset,\Xsigma)$
and any $\theta\in\Theta$,
\begin{eqnarray*}
&&
\bigl\llvert\Delta^{\theta}_{\chi,\chi'}\bigl\langle \chunk{z} {0} {n-1} \bigr\rangle
({f,h}) \bigr\rrvert
\\
&&\quad \leq\rho^{\lfloor n(\beta-\gamma^-)\rfloor} \bigl\{
\bigl(\chi\mathbf{L}^{\theta}\bigl\langle \chunk{z} {0} {n-1}
\bigr\rangle f \bigr) \bigl(\chi' \mathbf{L}^{\theta}\bigl\langle
\chunk{z} {0} {n-1} \bigr\rangle g\bigr) + \bigl(\chi'
\mathbf{L}^{\theta}\bigl\langle \chunk{z} {0} {n-1} \bigr\rangle
f\bigr) \bigl(\chi \mathbf{L}^{\theta}\bigl\langle \chunk{z} {0} {n-1}
\bigr\rangle g\bigr) \bigr\}
\\
&&\qquad{}
+ 2 \eta^{\lfloor n (\gamma^+ -\beta)\rfloor/2} \Biggl[ \prod_{i=0}^{n-1}
\bigl|\Pblock[\theta]{z_i}(\cdot, \Xset)\bigr|_{\infty}^2
\Biggr] |f|_{\infty} |h|_{\infty}.
\end{eqnarray*}
\end{prop}

\begin{pf}
Let $\eta>0$. According to~\ref{assumlikelihoodDrift}, there exists
a set $\Cset\subset\Yset$ such that (\ref{eqbound-eta-G}) and
(\ref{eqlower-bound}) hold. Denote
$\Csetprod\triangleq\Cset\times\Cset$ and for $z= \chunk{y}{0}{r-1}$,
set $\bar\varphi_{\Cset}^\theta\langle z \rangle=\varphi_{\Cset}^\theta
\langle z \rangle\otimes\varphi_{\Cset}^\theta\langle z
\rangle$ and $\probprod[\theta]{\Cset}{z}
\triangleq
\prob[\theta]{\Cset}{z} \otimes
\prob[\theta]{\Cset}{z}$ where
$\varphi_{\Cset}^\theta\langle z \rangle$ and $
\prob[\theta
]{\Cset}{z}$ are defined in Definition~\ref{defilocal-Doeblin}. For
any measurable nonnegative function $\fprod$ on
$(\Xsetprod,\Xsigmaprod)$, $\theta\in\Theta$ and $\xprod\in\Csetprod$,
%
\begin{eqnarray}
\label{eqminorization}
&&\bigl(\epsilon^-_{\Cset}(z)\bigr)^2
\bar\varphi_{\Cset}^\theta\langle z \rangle(\xprod) \probprod[
\theta]{\Cset} {z} (\one_{\Csetprod} \fprod)\nonumber\\[-8pt]\\[-8pt]
&&\qquad\leq\delta_{\xprod}
\Pblockprod[\theta]{z}(\one_{\Csetprod} \fprod)\leq\bigl(
\epsilon^{+}_{\Cset}(z)\bigr)^2 \bar
\varphi_{\Cset}^\theta\langle z \rangle(\xprod) \probprod[
\theta]{\Cset} {z}(\one_{\Csetprod} \fprod).\nonumber
\end{eqnarray}
Define the unnormalized kernel $\Pblockprod[\theta,0]{z}$ and
$\Pblockprod[\theta,1]{z}$
on $(\Xsetprod,\Xsigmaprod)$ as follows: for all $\xprod\in
\Xsetprod$
and $\bar{A} \in\Xsigmaprod$,
%
\begin{eqnarray}
\label{eqdefinition-Qprod0} \Pblockprod[\theta,0]{z}(\xprod,\bar
{A})&\triangleq&
\one_{\Csetprod}(\xprod) \bigl(\epsilon^-_\Cset(z)
\bigr)^2 \bar\varphi_{\Cset}^\theta\langle z \rangle(
\xprod) \probprod[\theta]{\Cset} {z} (\Csetprod\cap\bar{A} ),
\\
\label{eqdefinition-Qprod1} \Pblockprod[\theta,1]{z}(\xprod,\bar
{A})&\triangleq&
\Pblockprod[\theta]{z}(\xprod, \bar{A})- \Pblockprod[\theta
,0]{z}(\xprod,
\bar{A}).
\end{eqnarray}
Equation (\ref{eqminorization}) implies that, for all $\xprod\in
\Csetprod$, and any measurable nonnegative function~$\fprod$,
\[
0\leq\delta_{\xprod} \Pblockprod[\theta,1]{z}(\one_{\Csetprod}
\fprod) \leq r_\Cset(z) \delta_{\xprod}\Pblockprod[\theta]{z}(
\one_{\Csetprod
} \fprod),
\]
where
$r_\Cset(z)\triangleq1-(\epsilon^-_\Cset(z)/\epsilon^+_\Cset(z))^2$.
It then follows
%
\begin{eqnarray}
\label{eqP1Q}
&&\delta_{\xprod} \Pblockprod[\theta,1]{z}(\fprod)
\nonumber\\
&&\qquad= \one_{\Csetprod} (\xprod) \delta_{\xprod} \Pblockprod[
\theta,1]{z}( \one_{\Csetprod} \fprod) + \one_{\Csetprod}(\xprod)
\delta_{\xprod
}\Pblockprod[\theta,1]{z}(\one_{\Csetprod^c} \fprod) +
\one_{\Csetprod
^c}(\xprod) \delta_{\xprod}\Pblockprod[\theta,1]{z}(
\fprod)
\nonumber\\[-8pt]\\[-8pt]
&&\qquad\leq r_\Cset(z) \one_{\Csetprod} (\xprod)
\delta_{\xprod} \Pblockprod[\theta]{z}( \one_{\Csetprod} \fprod) +
\one_{\Csetprod}(\xprod) \delta_{\xprod} \Pblockprod[\theta]{z}(
\one_{\Csetprod^c} \fprod) + \one_{\Csetprod^c}(\xprod) \delta_{\xprod}
\Pblockprod[\theta]{z}(\fprod)
\nonumber\\
&&\qquad\leq\delta_{\xprod} \Pblockprod[\theta]{z}
\bigl(r_\Cset(z)^{\one
_{\Csetprod}(\xprod
) \one_{\Csetprod}} \fprod\bigr).\nonumber
\end{eqnarray}
Note that $\DDelta[\theta]{\chi,\chi'}{\chunk{z}{0}{n-1}}{f,h}$ may
be decomposed as
\[
\Delta^{\theta}_{\chi,\chi'}\bigl\langle \chunk{z} {0} {n-1} \bigr\rangle  ({f,h})= \sum
_{\chunk
{t}{0}{n-1} \in\{0,1\}^{n}} \Delta^{\theta,\chunk{t} {0} {n-1}}_{\chi,
\chi'}\bigl\langle \chunk{z} {0} {n-1} \bigr\rangle  ({f,h}),
\]
where
\[
\Delta^{\theta,\chunk{t} {0} {n-1}}_{\chi,\chi'}\bigl\langle \chunk{z} {0}
{n-1} \bigr\rangle  ({f,h})= \int\cdots\int\chi\otimes\chi'\bigl(\rmd
\bar{w}'_0\bigr) \Biggl( \prod
_{i=0}^{n-1} \Pblockprod[\theta,t_i]{z_i}(
\bar{w}_{i},\rmd\bar{w}_{i+1}) \Biggr) \Phi(
\bar{w}_{n})
\]
with $\Phi\triangleq f \otimes h-h \otimes f$. First assume that there
exists an index $i \in\{0,\ldots,n-1\}$ such that $t_i=0$. Then
\begin{eqnarray*}
\Delta^{\theta,\chunk{t} {0} {n-1}}_{\chi,\chi'}\bigl\langle \chunk{z} {0}
{n-1} \bigr\rangle  ({f,h}) &=&\chi\otimes\chi' \bigl(\Pblockprod[
\theta,t_0]{z_0} \cdots\Pblockprod[\theta,t_{i-1}]{z_{i-1}}
\bigl(\one_{\Csetprod}\times\bar\varphi_{\Cset}^\theta
\langle z_i \rangle\bigr) \bigr)
\\
&&{}\times\bigl(\epsilon^{-}_\Cset(z_i)
\bigr)^2 \probprod[\theta]{\Cset} {z_i} \bigl(
\one_{\Csetprod} \Pblockprod[\theta,t_{i+1}]{z_{i+1}}
\cdots\Pblockprod[\theta,t_{n-1}]{z_{n-1}} \Phi\bigr).
\end{eqnarray*}
By symmetry,
\[
\probprod[\theta]{\Cset} {z_i} \bigl( \one_{\Csetprod}
\Pblockprod[\theta,t_{i+1}]{z_{i+1}} \cdots\Pblockprod[
\theta,t_{n-1}]{z_{n-1}} \Phi\bigr) = 0,
\]
showing that $\DDelta[\theta,\chunk{t}{0}{n-1}]{\chi,\chi'}{\chunk
{z}{0}{n-1}}{f,h} = 0$ except if
for all $i \in\{0,\ldots,n-1\}$, $t_i=1$. Therefore,
\[
\Delta^{\theta}_{\chi,\chi'}\bigl\langle \chunk{z} {0} {n-1}
\bigr\rangle  ({f,h})= \chi \otimes\chi'
\bigl(\Pblockprod[\theta,1]{z_0} \cdots \Pblockprod[\theta,1]{z_{n-1}}
\Phi\bigr).
\]
This implies, using (\ref{eqP1Q}), that
%
\begin{eqnarray}\label{eqbound-G-Q}
&&
\bigl\llvert\Delta^{\theta}_{\chi,\chi'}\bigl\langle \chunk{z} {0} {n-1} \bigr\rangle
({f,h}) \bigr\rrvert\nonumber\\
&&\qquad\leq\chi\otimes\chi' \bigl(\Pblockprod[
\theta,1]{z_0} \cdots\Pblockprod[\theta,1]{z_{n-1}} |\Phi|
\bigr)
\nonumber\\[-8pt]\\[-8pt]
&&\qquad\leq\int\cdots\int\chi\otimes\chi'(\rmd\bar{w}_0)
\Biggl( \prod_{i=0}^{n-1}\Pblockprod[
\theta]{z_i}(\wprod_{i}, \rmd\wprod_{i+1})
\bigl( r_\Cset(z_i) \bigr)^{\one_{\Csetprod\times
\Csetprod
}(\wprod_{i},\wprod_{i+1})} \Biggr) \nonumber\\
&&\qquad\quad{}\times|\Phi|(\wprod_{n}).
\nonumber
\end{eqnarray}
Note that
%
\begin{equation}
\label{eqtool1} \prod_{i=0}^{n-1} \bigl(
r_\Cset(z_i) \bigr)^{\one_{\Csetprod
\times
\Csetprod}(\wprod_{i},\wprod_{i+1})} \leq
\varrho_\Cset^{\sum_{i=0}^{n-1} \one_{\Csetprod\times\Csetprod
}(\wprod
_{i},\wprod_{i+1})\one_{\Kset}(z_i)},
\end{equation}
where $\varrho_\Cset\triangleq\sup_{z \in\Kset} r_\Cset(z) < 1 $ under
\ref{assumlikelihoodDrift}. For any sequence $\chunk{z}{0}{n-1}$ such
that $n^{-1}\sum_{i=0}^{n-1} \one_{\Kset}(z_i) \geq(1-\gamma^-) $, we
have $\sum_{i=0}^{n-1} \one_{\Kset^c}(z_i) \leq n\gamma^- $, so that
\[
\sum_{i=0}^{n-1} \one_{\Kset^c}(z_i) \leq\bigl\lfloor n \gamma^-\bigr\rfloor.
\]
Moreover, we have
%
\begin{eqnarray}
\label{eqtool2}
&&
\sum_{i=0}^{n-1}
\one_{\Csetprod\times\Csetprod}(\wprod_{i},\wprod_{i+1})
\one_{\Kset}(z_i)
\nonumber\\
&&\qquad= \sum_{i=0}^{n-1}
\one_{\Csetprod\times\Csetprod}(\wprod_{i},\wprod_{i+1}) - \sum
_{i=0}^{n-1} \one_{\Csetprod\times\Csetprod}(
\wprod_{i},\wprod_{i+1}) \one_{\Kset^c}(z_i)
\nonumber\\[-8pt]\\[-8pt]
&&\qquad\geq N_{\Csetprod,n}(\chunk{\wprod} {0} {n}) - \sum
_{i=0}^{n-1} \one_{\Kset
^c}(z_i)
\nonumber\\
&&\qquad\geq N_{\Csetprod,n}(\chunk{\wprod} {0} {n}) - \bigl\lfloor n
\gamma^- \bigr\rfloor,\nonumber
\end{eqnarray}
where, for any set $\Asetprod\in\Xsigmaprod$, $N_{\Asetprod,n}(\chunk
{\wprod}{0}{n})= \sum_{i=0}^{n-1} \one_{\Asetprod\times\Asetprod
}(\wprod_{i},\wprod_{i+1})$. By combining (\ref{eqtool1}) and
(\ref{eqtool2}) and using that $\lfloor n \beta\rfloor- \lfloor n
\gamma^- \rfloor\geq\lfloor n(\beta-\gamma^-) \rfloor$, we therefore
obtain, for any $\beta\in\ocint{\gamma^-,1}$,
%
\begin{equation}
\label{eqbound-rho} \prod_{i=0}^{n-1}
\bigl( r_\Cset(z_i) \bigr)^{\one_{\Csetprod
\times
\Csetprod}(\wprod_{i},\wprod_{i+1})} \leq
\varrho_\Cset^{\lfloor n (\beta-\gamma^-) \rfloor} + \one\bigl\{
N_{\Csetprod,n}\bigl(\chunk{
\wprod} {0} {n}\bigr) < \lfloor n \beta\rfloor\bigr\}.
\end{equation}
For any sequence $\chunk{\wprod}{0}{n-1} \in\Xsetprod^{n}$ and any
$\Asetprod\in\Xsigmaprod$, denote
\[
M_{\Asetprod,n}\bigl(\chunk{\wprod} {0} {n-1}\bigr)\triangleq\sum
_{i=0}^{n-1} \one_{\Asetprod}(
\wprod_{i}).
\]
Using~\cite{doucfortmoulinespriouret2009}, Lemma 17, for any
sequence $\chunk{\wprod}{0}{n}$ satisfying $ N_{\Csetprod,n}(\chunk
{\wprod}{0}{n}) < \lfloor n \beta\rfloor$ which is equivalent to $
N_{\Csetprod,n}(\chunk{\wprod}{0}{n}) \leq\lfloor n \beta\rfloor
-1$, we have $M_{\Csetprod,n}(\chunk{\wprod}{0}{n-1}) \leq(\lfloor
n\beta\rfloor+n)/2$, so that
%
\begin{equation}
\label{eqlink-N-M} N_{\Csetprod,n}\bigl(\chunk{\wprod} {0} {n}\bigr) < \lfloor n
\beta
\rfloor\quad\Rightarrow\quad M_{\Csetprod^c,n}\bigl(\chunk{\wprod} {0} {n-1}\bigr) \geq
a_n \triangleq\frac{n-\lfloor n\beta\rfloor}{2}.
\end{equation}
In words, either the number of consecutive visits to the set $\Csetprod
$ at most $\lfloor n \beta\rfloor$, or the number of visits to the
complementary of the set $\Csetprod$ is larger than $a_n$. Plugging
(\ref{eqlink-N-M}) into (\ref{eqbound-rho}) and combining it with
(\ref{eqbound-G-Q}) yields
\begin{eqnarray*}
\bigl\llvert\Delta^{\theta}_{\chi,\chi'}\bigl\langle \chunk{z} {0} {n}
\bigr\rangle  ({f,h}) \bigr\rrvert&\leq&\varrho_\Cset^{\lfloor
n(\beta-\gamma^-) \rfloor} \chi \otimes \chi'
\bigl(\Pblockprod[\theta]{z_0} \cdots\Pblockprod[ \theta]{z_{n-1}}
|\Phi| \bigr)
\\
&&{}+ 2 |f|_{\infty} |h|_{\infty} \Gamma^{\theta
}_{\chi,\chi'}\bigl(
\chunk{z} {0} {n-1}\bigr),
\end{eqnarray*}
where
\[
\Gamma^\theta_{\chi,\chi'}\bigl(\chunk{z} {0} {n-1}\bigr)
\triangleq\int\cdots\int\chi\otimes\chi'(\rmd\wprod_0)
\prod_{i=0}^{n-1}\Pblockprod[
\theta]{z_i}(\wprod_{i}, \rmd\wprod_{i+1}) \one
\bigl\{ M_{\Csetprod
^c,n}\bigl(\chunk{\wprod} {0} {n-1}\bigr) \geq a_n \bigr\}.
\]
We finally have to bound this last term. First rewrite $\Gamma^\theta
_{\chi,\chi'}(\chunk{z}{0}{n-1})$ as follows:
\begin{eqnarray*}
\Gamma^\theta_{\chi,\chi'}(\chunk{z} {0} {n-1}) &=& \Biggl( \prod
_{i=0}^{n-1} \bigl|\Pblock[\theta]{z_i}(
\cdot,\Xset)\bigr|_{\infty}^2 \Biggr) \int\chi\otimes
\chi'(\rmd\wprod_0) \bigl(\eta^{\sum
_{i=0}^{n-1} \one_{\Csetprod^c}(\wprod_{i}) \one_{\Kset
}(z_i)} \bigr)
\\
&&{}\times\Biggl(\prod_{i=0}^{n-1}
\frac{\Pblockprod[\theta
]{z_i}(\wprod_{i},\rmd\wprod_{i+1}) }{\eta^{\one_{\Csetprod^c}(\wprod_{i})
\one_{\Kset
}(z_i)} |\Pblock[\theta]{z_i}(\cdot,\Xset)|_{\infty}^2 } \Biggr) \one
\bigl\{ M_{\Csetprod^c,n}\bigl(\chunk{\wprod} {0} {n-1}\bigr)
\geq a_n \bigr\}.
\end{eqnarray*}
Note that (\ref{eqdef-K}) implies that $\sum_{i=0}^{n-1} \one_{\Kset
}(z_i)\geq(n+\lfloor n \gamma^+ \rfloor)/2$. Then, for any $\gamma^+ >
\beta$, the inequality $M_{\Csetprod^c,n}(\chunk{\wprod}{0}{n-1})
\geq
a_n$ implies that
\[
\sum_{i=0}^{n-1} \one_{\Csetprod^c}(
\xprod_{i}) \one_{\Kset}(z_i) \geq\sum
_{i=0}^{n-1} \one_{\Csetprod^c}(
\xprod_{i})-\sum_{i=0}^{n-1}
\one_{\Kset
^c}(z_i)\geq\frac{\lfloor n \gamma^+\rfloor-\lfloor n \beta\rfloor
}{2} \geq\frac{\lfloor n(\gamma^+ - \beta)\rfloor}{2},
\]
showing that
\[
\bigl(\eta^{\sum_{i=0}^{n-1} \one_{\Csetprod^c}(\xprod_{i}) \one
_{\Kset
}(z_i)} \bigr) \one\bigl\{ M_{\Csetprod^c,n}\bigl(\chunk{\xprod}
{0} {n-1}\bigr) \geq a_n \bigr\} \leq\eta^{\lfloor n (\gamma^+ - \beta)
\rfloor/2}.
\]
The proof follows noting that, for any $\wprod=(w,w') \in\Xsetprod$
and $z \in\Yset^r$, (\ref{eqP-theta-vraisemblance}) and (\ref
{eqbound-eta-G}) imply
\[
\iint\frac{\Pblockprod[\theta]{z}(\wprod,\rmd\wprod_{i+1})
}{\eta^{\one
_{\Csetprod^c}(\wprod) \one_{\Kset}(z)} |\Pblock[\theta
]{z}(\cdot,\Xset)|_{\infty}^2 } =\frac{\Pblock[\theta
]{z}(w,\Xset)\Pblock[\theta
]{z}(w',\Xset) }{\eta^{\one_{\Csetprod^c}(\wprod) \one_{\Kset
}(z)} |\Pblock[\theta]{z}(\cdot,\Xset)|_{\infty}^2 }\leq1.
\]
\upqed\end{pf}
%
%
\begin{lem}\label{lemboundFilter}
Let $(U_k)_{k \in\zset}, (V_k)_{k \in\zset}, (W_k)_{k \in\zset}$ be
stationary sequences such that
\[
\mathbb E\bigl[{\ln^+ U_0}\bigr]<\infty,\qquad \mathbb E\bigl[{\ln^+
V_0}\bigr]<\infty,\qquad \mathbb E\bigl[{\ln^+ W_0}\bigr]<\infty.
\]
Then, for all $\eta,\rho$ in $(0,1)$ such that $-\ln\eta>\PE{\ln^+
V_0}$, there exists a $\PP{}$-a.s. finite random variable $D$ and a
constant $\varrho\in(0,1)$ such that for all $k \geq1, m\geq0$,
\[
\rho^{k+m}+\eta^{k+m} W_{-m} \Biggl(\prod
_{i=-m}^{k-1} V_i \Biggr) U_k
\leq\varrho^{k+m} D,\qquad \PP{}\mbox{-a.s.}
\]
\end{lem}
\begin{pf} Let $\alpha\in(0,1)$ such that $\PE{\ln^+ V_0}<-\ln
\alpha
<-\ln\eta$, and let $\tilde{\alpha}>0$ such that $(\eta/\alpha)
\vee
\rho<\tilde{\alpha}<1$. Then
\begin{eqnarray*}
&&\rho^{k+m}+\eta^{k+m} W_{-m} \Biggl(\prod
_{i=-m}^{k-1} V_i \Biggr) U_k
\\
&&\qquad= \Biggl[ \biggl(\frac{\rho}{\tilde{\alpha}} \biggr)^{k+m}\tilde
{\alpha
}^m+ \biggl( \frac{\eta}{\alpha\tilde{\alpha}} \biggr)^{k+m} \bigl(
\tilde
{\alpha}^m W_{-m} \bigr) \Biggl( \prod
_{i=-m}^{k-1} (V_i\alpha) \Biggr) \bigl(
\tilde{\alpha}^k U_k \bigr) \Biggr]
\\
&&\qquad\leq\biggl( \frac{\rho}{\tilde{\alpha}} \vee\frac{\eta}{\alpha
\tilde
{\alpha}} \biggr)^{k+m}
D
\end{eqnarray*}
with
\begin{eqnarray*}
D &\triangleq& 1
+ \Bigl(\sup_{m\geq0} \tilde{\alpha}^m W_{-m}
\Bigr) \Biggl(\sup_{m\geq
0} \prod_{i=-m}^0
(V_i\alpha) \Biggr) \Biggl(\sup_{k\geq1} \prod
_{i=1}^{k-1} (V_i\alpha) \Biggr) \Bigl(
\sup_{k\geq1} \tilde{\alpha}^k U_k \Bigr).
\end{eqnarray*}
We now show that $D$ is $\PP{}$-a.s. finite. First note that combining
the bound $\PE{\ln^+ U_0<\infty}$ with Lemma~\ref{lemstationary}
(stated and proved below), we obtain that the random variable $\sup
_{k\geq1} \tilde{\alpha}^k U_k$ is $\PP{}$-a.s. finite; in the same
way, $\sup_{m\geq0} \tilde{\alpha}^m W_{-m}$ is $\PP{}$-a.s. finite.
Moreover, since $\PE{\ln^+ V_0}<\infty$, Birkoff's ergodic theorem
ensures that
\[
\frac{1}{k-1} \sum_{i=1}^{k-1} \ln^+
V_i \rightarrow_{k \to\infty
} \PE{\ln^+ V_0} <-{\ln
\alpha}, \qquad\PP{}\mbox{-a.s.}
\]
By taking the exponential function in the previous limit, we obtain that
\[
\prod_{i=1}^{k-1} (V_i\alpha)
\leq\exp\Biggl\{(k-1) \Biggl(\frac{1}{k-1} \sum
_{i=1}^{k-1} \ln^+ V_i + \ln\alpha
\Biggr) \Biggr\}\rightarrow_{k
\to\infty}0,\qquad \PP{}\mbox{-a.s.}
\]
so that $\sup_{k \geq1}\prod_{i=1}^{k-1} (V_i\alpha)$ is $\PP{}$-a.s.
finite. Following the same arguments,
\[
\sup_{m\geq0} \prod_{i=-m}^0
(V_i\alpha)
\]
is $\PP{}$-a.s. finite. Finally $D$ is $\PP{}$-a.s. finite.
The proof is complete.
\end{pf}
%
\begin{lem}\label{lemstationary}
Let $\{Z_k\}_{k \in\zset}$ be a sequence of nonnegative random
variables on a probability space $(\Omega,\mathcal{A},\PP{})$ having
the same marginal distribution, that is, for any $k \in\zset$ and any
measurable nonnegative function $f$, $\PE{f(Z_k)}=\PE{f(Z_0)}$.
\renewcommand\thelonglist{(\roman{longlist})}
\renewcommand\labellonglist{\thelonglist}
\begin{longlist}
\item Assume that $\PE{(\ln Z_0)^+}<\infty$. Then, for all $\beta\in
(0,1)$, $\sup_{k\geq0} \beta^k Z_k<\infty$, $\PP{}$-a.s.
\item Assume that $\PE{|\ln Z_0|}<\infty$. Then, for all $\beta\in
(0,1)$, $\sup_{k\in\zset} \beta^{|k|} Z_k<\infty$ and $\inf_{k\in
\zset
} \beta^{-|k|} Z_k>0$, $\PP{}$-a.s.
\end{longlist}
\end{lem}
\begin{pf}
Let $\beta\in(0,1)$. Since
\[
\mathbb P\bigl[{ \beta^k Z_k>1 }\bigr] = \mathbb P\bigl[{\ln Z_k / (-{\ln
\beta}) \geq k}\bigr] = \mathbb P\bigl[\ln Z_0 / (-{\ln\beta}) \geq k\bigr],
\]
it follows that
\[
\sum_{k = 0}^\infty\mathbb P\bigl[{\beta^k
Z_k>1 }\bigr] = \sum_{k = 0}^\infty\mathbb P\bigl[{
\ln Z_0 / (-{\ln\beta}) \geq k }\bigr]\leq\mathbb E\bigl[{(\ln Z_0)^+}\bigr]/ (-{\ln
\beta}) < \infty.
\]
The proof of (i) is completed by using the Borel--Cantelli lemma. Now,
(ii) can be easily derived by noting that if $\PE{|{\ln Z_0}|}<\infty$,
then one may use twice (i), first by replacing $Z_k$ by $Z_{-k}$ and
then by replacing $Z_k$ by $1/Z_k$.~%
\end{pf}
%
\begin{prop}
\label{propforgetting-initial-condition-likelihood} Assume
\ref{assumlikelihoodDrift} and~\ref{assummajo-g}. There exist a constant
$\kappa\in(0,1)$, an integer-valued random variable $K$ satisfying
$\PP[Y]{K < \infty}= 1$ such that, for any initial distributions $\chi,
\chi' \in\mathcal{M}(\Dset,r)$ [where $\mathcal{M}(\Dset,r)$ is defined
in (\ref{eqmort-de-rire})],
%
\begin{eqnarray}
\label{itemforgetting-initial-condition-loglikelihood}
\sup_{\theta\in
\Theta}
\sup_{k \geq K}\sup_{m \geq0} \kappa^{-(m+k)} \bigl\llvert\ln
 p_{\chi}^{\theta}\bigl(Z_k|\chunk{Z} {-m} {k-1}\bigr)
-\ln p_{\chi'}^{\theta}\bigl(Z_k|\chunk{Z} {-m} {k-1}\bigr) \bigr
\rrvert&<&\infty,\nonumber\\[-8pt]\\[-8pt]
&&\eqntext{\PP{}\mbox{-a.s.},}
\\
\label{itemforgetting-initial-condition-loglikelihood-vers-moins-infini}
\sup_{\theta\in\Theta} \sup_{k \geq K}\sup_{m \geq0}
\kappa^{-(m+k)} \bigl\llvert\ln p_{\chi}^{\theta}\bigl(Z_k|
\chunk{Z} {-m} {k-1}\bigr) -\ln p_{\chi}^{\theta}\bigl(Z_k|\chunk{Z} {-m-1}
{k-1}\bigr) \bigr\rrvert&<&\infty,\nonumber\\[-8pt]\\[-8pt]
&&\eqntext{\PP{}\mbox{-a.s.},}
\\
\label{itemforgetting-initial-condition-loglikelihood-vers-moins-infini-k=0}
\sup_{\theta\in\Theta}\sup_{m \geq0} \kappa^{-m} \bigl\llvert
\ln p_{\chi}^{\theta}\bigl(Z_0|\chunk{Z} {-m} {-1}\bigr) -\ln p_{\chi}^{\theta}
\bigl(Z_0|\chunk{Z} {-m-1} {-1}\bigr) \bigr\rrvert&<&\infty,\nonumber\\[-8pt]\\[-8pt]
&&\eqntext{\PP{}
\mbox{-a.s.}}
\end{eqnarray}
\end{prop}

\begin{pf}
\textit{Proof of}
(\ref{itemforgetting-initial-condition-loglikelihood}).
It follows from (\ref{eqvrais-en-fonction-de-P}) that, for any integer
$(m,k) \in\nset$
and any sequence $\chunk{z}{-m}{k}$,
\[
p_{\chi}^{\theta}\bigl(z_k|\chunk{z} {-m} {k-1}\bigr)=
\frac{\chi\Pblock
[\theta]{\chunk{z}{-m}{k-1}} (\Pblock[\theta]{z_k}\one_\Xset
)}{\chi
\Pblock[\theta]{\chunk{z}{-m}{k-1}} (\one_\Xset)}.
\]
Since, for any $a,b > 0$, $\ln(a) - \ln(b) \leq(a-b)/b$, definition
(\ref{eqdefinition-Delta}) implies that
%
\begin{eqnarray}\label{eqbound-ln}
&&
\ln p_{\chi}^{\theta}\bigl(z_k|\chunk{z} {-m} {k-1}\bigr)
- \ln p_{\chi'}^{\theta}\bigl(z_k|\chunk{z} {-m}
{k-1}\bigr)\nonumber\\[-8pt]\\[-8pt]
&&\qquad
\leq\frac{\DDelta[\theta]{\chi,\chi'}{\chunk{z}{-m}{k-1}}{\Pblock
[\theta
]{z_k}\one_\Xset,\one_\Xset}}{\chi\Pblock[\theta]{\chunk
{z}{-m}{k-1}}(\one_\Xset)\times\chi' \Pblock[\theta]{\chunk
{z}{-m}{k-1}}(\Pblock
[\theta
]{z_k}\one_\Xset)}.\nonumber
\end{eqnarray}
Let $0 \leq\gamma^- < \gamma^+ \leq1$.
By Proposition~\ref{propmajo-difference-numerateur}, for any $\eta>
0$ and $\beta\in\ooint{\gamma^-,\gamma^+}$ there exists $\varrho
\in
(0,1)$ such that, for any sequence $\chunk{z}{-m}{k-1}$
satisfying
%
\begin{equation}
\label{eqcondition-sur-visite-de-z} (m+k)^{-1} \sum
_{i=-m}^{k-1} \one_{\Kset}(z_i)
\geq\max\bigl(1-\gamma^-,\bigl(1+\gamma^+\bigr)/2\bigr),
\end{equation}
we have
%
\begin{eqnarray}
\label{eqpremiere-ineg}
&&\frac{\DDelta[\theta]{\chi,\chi'}{\chunk
{z}{-m}{k-1}}{\Pblock
[\theta]{z_k}\one_\Xset,\one_\Xset}}{\chi\Pblock[\theta]{\chunk
{z}{-m}{k-1}}(\one_\Xset)\times\chi' \Pblock[\theta]{\chunk
{z}{-m}{k-1}}(\one_\Xset)}
\nonumber\\
&&\qquad
\leq\varrho^{a(m+k)} \biggl[ 1 + \frac{\chi\Pblock[\theta]{\chunk
{z}{-m}{k-1}}(\Pblock[\theta]{z_k}\one_\Xset) \times\chi'
\Pblock
[\theta]{\chunk{z}{-m}{k-1}}(\one_\Xset)} {\chi\Pblock[\theta
]{\chunk
{z}{-m}{k-1}}(\one_\Xset) \times\chi'\Pblock[\theta]{\chunk
{z}{-m}{k-1}}(\Pblock[\theta]{z_k}\one_\Xset)} \biggr]\\
&&\qquad\quad{} +
2 \eta^{b(m+k)} C_{m,k},\nonumber
\end{eqnarray}
where $a(n)=\lfloor n (\beta-\gamma^-)\rfloor$, $b(n)=\lfloor
n(\gamma^+-\beta)\rfloor/2$ and
%
\begin{equation}
\label{eqdefinition-Cm,k} C_{m,k} \triangleq
\frac{\prod_{i=-m}^{k-1} |\Pblock[\theta
]{z_i}(\cdot,\Xset)|_{\infty}^2 } {\chi\Pblock[\theta]{\chunk
{z}{-m}{k-1}}(\one_\Xset) \times\chi'\Pblock[\theta]{\chunk
{z}{-m}{k-1}}(\Pblock
[\theta]{z_k}\one_\Xset)} \bigl|\Pblock[\theta]{z_k}(\cdot,\Xset
)\bigr|_{\infty}.
\end{equation}
Moreover, by (\ref{eqdefinition-Delta}),
\begin{eqnarray*}
&&
\frac{\chi\Pblock[\theta]{\chunk{z}{-m}{k-1}}(\Pblock[\theta
]{z_k}\one_\Xset) \times\chi' \Pblock[\theta]{\chunk{z}{-m}{k-1}}(\one
_\Xset
)} {\chi\Pblock[\theta]{\chunk{z}{-m}{k-1}}(\one_\Xset) \times
\chi'\Pblock[\theta]{\chunk{z}{-m}{k-1}}(\Pblock[\theta]{z_k}\one_\Xset)}
\\
&&\qquad
= \frac{\DDelta[\theta]{\chi,\chi'}{\chunk{z}{-m}{k-1}}{\Pblock
[\theta
]{z_k} \one_\Xset,\one_{\Xset}}} {\chi\Pblock[\theta]{\chunk
{z}{-m}{k-1}}(\one_\Xset) \times\chi'\Pblock[\theta]{\chunk
{z}{-m}{k-1}}(\Pblock[\theta]{z_k}\one_\Xset)} + 1.
\end{eqnarray*}
Plugging this identity into (\ref{eqpremiere-ineg}) and then using
(\ref{eqbound-ln}) yields
%
\begin{eqnarray}\label{eqinegalite-log-dens-bloc}
&&
\ln p_{\chi}^{\theta}\bigl(z_k|\chunk{z} {-m} {k-1}\bigr) - \ln
p_{\chi'}^{\theta}\bigl(z_k|\chunk{z} {-m} {k-1}\bigr) \nonumber\\[-8pt]\\[-8pt]
&&\qquad\leq 2
\bigl(1-\varrho^{a(m+k)}\bigr)^{-1}
\bigl[ \varrho^{a(m+k)} + \eta^{b(m+k)} C_{m,k}
\bigr].
\nonumber
\end{eqnarray}
For any sequence $\chunk{z}{-m}{k-1}$, we have
%
\begin{eqnarray}\label{eqmino-denom}
\chi\mathbf{L}^{\theta}\bigl\langle \chunk{z} {-m} {k-1} \bigr\rangle (\one_\Xset) &\geq&\chi(
\Dset) \prod_{i=-m}^{k-1} \Bigl\{
\inf_{x \in\Dset} \Pblock[\theta]{z_i}(x,\Dset) \Bigr\},
\nonumber\\[-8pt]\\[-8pt]
\chi'\mathbf{L}^{\theta}\bigl\langle \chunk{z} {-m} {k-1} \bigr\rangle \bigl(\Pblock[\theta
]{z_k}\one_\Xset\bigr) &\geq&\chi'(\Dset) \prod
_{i=-m}^{k} \Bigl\{ \inf_{x \in\Dset}
\Pblock[\theta]{z_i}(x,\Dset) \Bigr\}.
\nonumber
\end{eqnarray}
Exchanging $\chi$ and $\chi'$ in (\ref{eqinegalite-log-dens-bloc})
allows us to obtain an upper bound for\break $|{\ln
p_{\chi}^{\theta}}(z_k|\chunk{z}{-m}{k-1}) -
\ln\cdens[\chi'][\theta]{z_k}{\chunk {z}{-m}{k-1}}|$. More precisely,
for any sequence $\chunk{z}{-m}{k-1}$ satisfying
(\ref{eqcondition-sur-visite-de-z}), we have
%
\begin{eqnarray}\label{eqbound-diff-ln-vers-moins-infini}
&&
\sup_{\theta\in\Theta} \bigl\llvert\ln p_{\chi}^{\theta}\bigl(z_k| \chunk{z}
{-m} {k-1}\bigr) - \ln p_{\chi'}^{\theta}\bigl(z_k|\chunk{z} {-m} {k-1}\bigr)
\bigr\rrvert\nonumber\\
&&\qquad\leq 2\bigl(1- \varrho^{a(m+k)}\bigr)^{-1}
\\
&&\qquad\quad{}\times\Biggl\{ \varrho^{a(m+k)} + \frac{\eta^{b(m+k)}}{\chi(\Dset)
\chi'(\Dset)} \Biggl[\prod
_{j=-m}^{k-1} (D_{z_j}
)^2 \Biggr] D_{z_k} \Biggr\},
\nonumber
\end{eqnarray}
where, for $z\in\Yset^r $,
%
\begin{equation}
\label{eqdefinition-D-z} D_z= \frac{\sup_{\theta\in\Theta}
|\Pblock[\theta]{z}(\cdot,\Xset)|_{\infty}} {\inf_{\theta\in\Theta}
\inf_{x \in\Dset}\Pblock[\theta]{z}(x,\Dset)}.
\end{equation}
Assume that $\PE{\ln^+ (D_{Z_0}) }<\infty$, and set $\eta$ small enough
so that $\PE{\ln^+ (D_{Z_0})} \leq- \ln\eta$. By Lemma \ref
{lemboundFilter}, there exists a $\PP{}$-a.s. finite random variable
$C$, and a constant $\kappa\in(0,1)$ such that, for
all $k \geq1$, $m\geq0$,
\[
\frac{2}{1-\varrho^{a(m+k)}} \Biggl\{ \varrho^{a(m+k)} + \frac{\eta
^{b(m+k)}}{\chi(\Dset) \chi'(\Dset)} \Biggl[
\prod_{j=-m}^{k-1} (D_{z_j}
)^2 \Biggr] D_{z_k} \Biggr\} \leq C \kappa^{k+m},\qquad
\PP{}\mbox{-a.s.}
\]
It remains to show that $\PE{\ln^+(D_{Z_0}) }<\infty$. Since for any
$a,b>0$, $\ln^+(a/b)\leq\ln^+ (a) + \ln^- (b)$,
%
\begin{equation}
\label{eqmajo-ln-D} \ln^+ (D_{z} ) \leq\ln^+ \Bigl(
\sup_{\theta\in\Theta} \bigl|\Pblock[\theta]{z}(\cdot,\Xset)\bigr|_{\infty} \Bigr
) +
\ln^- \Bigl(\inf_{\theta\in\Theta}\inf_{x \in\Dset} \Pblock[\theta
]{z}(x,\Dset)
\Bigr).
\end{equation}
Since,\vspace*{1pt} for any $z = \chunk{y}{0}{r-1} \in\Yset^r$, $\sup_{\theta
\in
\Theta} |\Pblock[\theta]{z}(\cdot,\Xset)|_{\infty} \leq\prod
_{i=0}^{r-1} \sup_{\theta\in\Theta} |\g[\theta](\cdot,y_i)|_{\infty}$,
\ref{assumlikelihoodDrift}\ref{itemcondition-minoration} and \ref
{assummajo-g} imply that $\PE{\ln^+ (D_{Z_0}) }<\infty$.
Finally, according to (\ref{eqbound-diff-ln-vers-moins-infini}),
\[
\sup_{\theta\in\Theta}\bigl\llvert\ln p_{\chi}^{\theta}\bigl(Z_k|
\chunk{Z} {-m} {k-1}\bigr) -\ln p_{\chi'}^{\theta}\bigl(Z_k|\chunk{Z} {-m} {k-1}\bigr)
\bigr\rrvert\leq C
\kappa^{m+k},\qquad \PP{}\mbox{-a.s.},
\]
provided that
%
\begin{equation}
\label{eqconditionVisiteC} (m+k)^{-1} \sum
_{j=-m}^{k-1} \one_{\Kset}(Z_j)
\geq\max\bigl(1-\gamma^-,\bigl(1+\gamma^+\bigr)/2\bigr),\qquad \PP{}\mbox{-a.s.}
\end{equation}
It thus remains to show the existence of a $\PP{}$-a.s. finite random
variable $K$ such that for any $k\geq K$ and any $m \geq0$, (\ref
{eqconditionVisiteC}) holds $\PP{}$-a.s.
Under~\ref{assumlikelihoodDrift}\ref{itemcondition-L-K}, $1 - \PP
{Z_0 \in\Kset} < 2 \PP{Z_0 \in\Kset} -1$. Then, choose
$\tilde\gamma^-$, $\gamma^-$, $\gamma^+$ and $\tilde\gamma^+$
such that
%
\begin{equation}
\label{eqcondition-tilde-gamma-vs-gamma} 1 - \PP{Z_0 \in\Kset} <
\tilde\gamma^-<\gamma^- < \gamma^+ < \tilde\gamma^+<2 \PP{Z_0 \in
\Kset} -1.
\end{equation}
By construction $ (1+\tilde{\gamma}^+)/2 < \PP[Y]{Z_0 \in\Kset}$
and $1
- \tilde{\gamma}^- < \PP{Z_0 \in\Kset}$.
Since $(Z_k)_{k \in\zset}$ is stationary and ergodic, the Birkhoff
ergodic theorem ensures that there exists a $\PP{}$-a.s. finite random
variable $B$ such that for any $k \geq B$ and $m \geq B$, $\PP{}$-a.s.,
%
\begin{eqnarray}
\max\biggl(1-\tilde{\gamma}{}^-,\frac{1+\tilde{\gamma}^+}{2} \biggr) &<& k^{-1}
\sum_{i=0}^{k-1} \one_\Kset(Z_i),
\\
\label{eqbirkhoff-somme-m-one}
\max\biggl(1-\tilde{\gamma}^-,\frac{1+\tilde{\gamma}^+}{2} \biggr) &<& m^{-1}
\sum_{i=-m}^{-1} \one_{\Kset}(Z_i).
\end{eqnarray}
Set $K^+\triangleq B (1+\gamma^+)/(\tilde{\gamma}^+ - \gamma^+)$.
If $m\geq
B$ and $k \geq K^+$, then using that $K^+ \geq B$, $\PP{}$-a.s.,
\[
\frac{\sum_{i=-m}^{k-1} \one_\Kset(Z_i)}{k+m}> \frac{k(1+\tilde
{\gamma
}^+)/2+m(1+\tilde{\gamma}^+)/2}{k+m}=\bigl(1+\tilde{\gamma}^+\bigr)/2>
\bigl(1+\gamma^+\bigr)/2.
\]
Now, if $0 \leq m< B$ and $k \geq K^+$,
\begin{eqnarray*}
\frac{\sum_{i=-m}^{k-1} \one_\Kset(Z_i)}{k+m}&\geq&\frac{\sum_{i=0}^{k-1}
\one_\Kset(Z_i)}{k+m}
>\frac{k(1+\tilde{\gamma}^+)/2}{k+m} \\
&>& \frac{ K^+(1+\tilde{\gamma
}^+)/2}{ K^++B}=\bigl(1+\gamma^+\bigr)/2.
\end{eqnarray*}
Similarly, setting $K^- \triangleq B(1-\gamma^-)/(\tilde{\gamma}^- -
\gamma^-)$, we obtain, for all $m \geq0$ and all $k \geq K^-$ that,
$\PP{}$-a.s.,
\[
\frac{\sum_{i=-m}^{k-1} \one_\Kset(Z_i)}{k+m} \geq1 - \gamma^-.
\]
The proof of (\ref{itemforgetting-initial-condition-loglikelihood}) is
now completed by setting $K = K^+ \vee K^-$.\vspace*{9pt}

\textit{Proof of} (\ref
{itemforgetting-initial-condition-loglikelihood-vers-moins-infini}).
Note that
\[
p_{\chi}^{\theta}\bigl(z_k|\chunk{z} {-m-1} {k-1}\bigr)=
p_{\chi'}^{\theta}\bigl(z_k|\chunk{z} {-m} {k-1}\bigr)
\]
with $\chi'(\Aset)=\chi(\Pblock[\theta]{z_{-m-1}} \one_\Aset
)/\chi
(\Pblock[\theta]{z_{-m-1}} \one_\Xset)$. Since
\[
\frac{1}{\chi'(\Dset)}=\frac{\chi(\Pblock[\theta]{z_{-m-1}} \one_\Xset
)}{\chi(\Pblock[\theta]{z_{-m-1}} \one_\Dset)} \leq\frac
{D_{z_{-m-1}}}{\chi(\Dset)},
\]
where $D_z$ is defined in (\ref{eqdefinition-D-z}), (\ref
{eqbound-diff-ln-vers-moins-infini}) writes
\begin{eqnarray*}
&&
\sup_{\theta\in\Theta} \bigl\llvert\ln p_{\chi}^{\theta}\bigl(z_k|
\chunk{z} {-m} {k-1}\bigr) - \ln p_{\chi}^{\theta}\bigl(z_k|\chunk{z}
{-m-1} {k-1}\bigr) \bigr\rrvert\\
&&\qquad\leq 2\bigl(1-\varrho^{a(m+k)}
\bigr)^{-1}
\\
&&\qquad\quad{}\times\Biggl[ \varrho^{a(m+k)} + \frac{\eta^{b(m+k)}}{[\chi(\Dset
)]^2} D_{z_{-m-1}}
\prod_{j=-m}^{k-1} (D_{z_j}
)^2 D_{z_k} \Biggr].
\end{eqnarray*}
And the rest of the proof of (\ref
{itemforgetting-initial-condition-loglikelihood-vers-moins-infini})
follows the same lines
as (\ref{itemforgetting-initial-condition-loglikelihood}) and is
omitted for brevity.\vspace*{9pt}

\textit{Proof of} (\ref
{itemforgetting-initial-condition-loglikelihood-vers-moins-infini-k=0}).
Noting that, when $k=0$, equation (\ref{eqconditionVisiteC})
follows immediately from~(\ref{eqbirkhoff-somme-m-one}),
the proof of (\ref
{itemforgetting-initial-condition-loglikelihood-vers-moins-infini-k=0})
follows the same lines as the proof of (\ref
{itemforgetting-initial-condition-loglikelihood-vers-moins-infini}) and
is omitted for brevity.
\end{pf}
%
\begin{cor}[(Corollary of Proposition~\ref{propforgetting-initial-condition-likelihood})]
Assume~\ref{assumlikelihoodDrift} and~\ref{assummajo-g}.
For any $\theta\in\Theta$, there exists a measurable function $\pi
^\theta_Z\dvtx  \Zset^{\zset^{-}} \to\rset$ such that for any probability
measure $\chi$ satisfying $\chi(\Dset) \in\mathcal{M}(\Dset,r)$\vadjust{\goodbreak}
[where $\mathcal{M}(\Dset,r)$ is defined in (\ref{eqmort-de-rire})],
%
\begin{equation}
\label{eqexistence-approx-station} \mathbb P_{Y}\Bigl[{\lim_{m \to\infty}
p_{\chi}^{\theta}\bigl(Z_0|\chunk{Z} {-m} {-1}\bigr)=\pi^\theta_Z\bigl(
\chunk{Z} {-\infty} {0}\bigr)}\Bigr]=1.
\end{equation}
\end{cor}
In the sequel, we denote $\dens[\ ][\theta]{Z_0|\chunk{Z}{-\infty}{-1}}
\triangleq\pi^\theta_Z(\chunk{Z}{-\infty}{0})$ and for $n
\geq0$,\break
$p^{\theta}(Z_{0}^{n}|\chunk{Z}{-\infty}{-1}) \triangleq
\prod_{i=0}^n \pi^\theta_Z(\chunk{Z}{-\infty}{i})$.

\subsubsection{Consistency of the block MLE}
%
\begin{prop} \label{proplimiteRelativeEntropyRateZ}
Assume~\ref{assumlikelihoodDrift} and~\ref{assummajo-g}. Then:

\renewcommand\thelonglist{(\roman{longlist})}
\renewcommand\labellonglist{\thelonglist}
\begin{longlist}
\item\label{itemexistence-approx-stat-bloc} For any $\theta\in
\Theta$,
%
\begin{equation}
\label{eqesperance-log-stat-bloc} \mathbb E\bigl[\bigl|\ln p^{\theta}\bigl(Z_0|
\chunk{Z} {-\infty} {-1}\bigr)\bigr|\bigr]<\infty.
\end{equation}
\item\label{itemconv-unif-entropy-stat-bloc}For any probability
measure $\chi\in\mathcal{M}(\Dset,r)$ [where $\mathcal{M}(\Dset,r)$
is defined in~(\ref{eqmort-de-rire})],
\[
\limsup_{n \to\infty} \sup_{\theta\in\Theta} \bigl\llvert n^{-1} \ln
p_{\chi}^{\theta}\bigl(Z_{0}^{n-1}\bigr)- n^{-1}
\ln p^{\theta}\bigl(Z_{0}^{n-1}|\chunk{Z} {-\infty} {-1}\bigr) \bigr\rrvert=0,
\qquad\PP{}
\mbox{-a.s.}
\]
\item\label{itemconv-ponct-entropy-bloc}For any $\theta\in\Theta$,
and for any probability measure $\chi\in\mathcal{M}(\Dset,r)$,
\[
\lim_{n \to\infty} n^{-1} \ln p_{\chi}^{\theta}\bigl(Z_{0}^{n-1}\bigr)
=\mathbb E\bigl[\ln p^{\theta}\bigl(Z_0|\chunk{Z} {-\infty} {-1}\bigr)\bigr],\qquad
\PP{}\mbox{-a.s.}
\]
\end{longlist}
\end{prop}
\begin{pf}
\textit{Proof of}~\ref{itemexistence-approx-stat-bloc}.
It follows from (\ref{eqexistence-approx-station}) that, $\PP{}$-a.s.,
%
\begin{equation}
\label{eqmajo-cond-dens-par-G}\qquad
p^{\theta}\bigl(Z_0|\chunk{Z} {-
\infty} {-1}\bigr) =\lim_{m \to\infty}  p_{\chi}^{\theta}\bigl(Z_0|\chunk
{Z} {-m} {-1}\bigr)\leq\bigl|\Pblock[\theta]{Z_0}(\cdot,\Xset)\bigr|_{\infty}
\leq\prod_{i=0}^{r-1} \bigl|\g[\theta](
\cdot,Y_i)\bigr|_{\infty}.
\end{equation}
Then,~\ref{assummajo-g} shows that
\[
\mathbb E\bigl[\ln^+  p^{\theta}\bigl(Z_0|\chunk{Z} {-\infty} {-1}\bigr)\bigr]
\leq\mathbb E\bigl[
\ln^+ \bigl|\Pblock[\theta]{Z_0}(\cdot,\Xset)\bigr|_{\infty} \bigr]<\infty.
\]
We now show that $\PE{\ln^- \dens[\ ][\theta]{Z_{0}|\chunk
{Z}{-\infty
}{-1}}} < \infty$ by establishing that\break
$\mathbb E[\ln p^{\theta}(Z_{0}|\chunk{Z}{-\infty}{-1})] >-\infty$.
For that purpose, introduce the sequence
\[
L_m^\theta\triangleq m^{-1} \sum _{\ell=1}^{m} \bigl[\ln^+
\bigl|\Pblock[ \theta]{Z_0}(\cdot,\Xset)\bigr|_{\infty} -\ln
p_{\chi}^{\theta}\bigl(Z_{0}|\chunk{Z} {-\ell} {-1}\bigr) \bigr].
\]
By (\ref{eqmajo-cond-dens-par-G}), the sequence $(L_m^\theta)_{m
\geq
0}$ is nonnegative and the Fatou lemma implies that
%
\begin{equation}
\label{eqfatou} \liminf_{m \to\infty} \mathbb E\bigl[{L_m^\theta}\bigr]
\geq\mathbb E\Bigl[{\liminf_{m \to
\infty} L_m^\theta}\Bigr].
\end{equation}
By definition,
%
\begin{eqnarray}
\label{eqliminfPE} \liminf_{m \to\infty} \mathbb E\bigl[{L_m^\theta}\bigr]
&=&\mathbb E\bigl[{\ln^+ \bigl|\Pblock[\theta]{Z_0}(\cdot,\Xset)\bigr|_{\infty} }\bigr]
\nonumber\\[-8pt]\\[-8pt]
&&{}-\limsup_{m \to\infty} m^{-1} \sum_{\ell=1}^{m}
\mathbb E\bigl[\ln p_{\chi}^{\theta}\bigl(Z_{0}|\chunk{Z} {-\ell}
{-1}\bigr)\bigr]
\nonumber
\end{eqnarray}
and
%
\begin{eqnarray}
\label{eqPEliminf} \mathbb E\Bigl[{\liminf_{m \to\infty} L_m^\theta}\Bigr] &=&
\mathbb E\bigl[{\ln^+ \bigl|\Pblock[\theta]{Z_0}(\cdot,\Xset)\bigr|_{\infty} }\bigr]
\nonumber\\[-8pt]\\[-8pt]
&&{}- \mathbb E\Biggl[\limsup_{m \to\infty} m^{-1} \sum
_{\ell= 1}^{m} \ln p_{\chi}^{\theta}\bigl(Z_{0}|
\chunk{Z} {-\ell} {-1}\bigr)\Biggr].
\nonumber
\end{eqnarray}
Since $(Y_k)_{k \in\zset}$ is stationary, for any \mbox{$\ell\in\nset,\hspace*{-0.1pt}
\PE{ \ln\dens[\chi][\theta]{Z_{0}|\chunk{Z}{-\ell}{-1}}}\hspace*{-0.1pt}= \hspace*{-0.1pt}\PE{
\ln\dens[\chi][\theta]{Z_{\ell}|\chunk{Z}{0}{\ell-1}}}$}
showing that
%
\begin{equation}
\label{eqidentity-1} m^{-1} \sum_{\ell=1}^{m}
\mathbb E\bigl[{\ln p_{\chi}^{\theta}\bigl(Z_{0}|\chunk{Z} {-\ell} {-1}\bigr)}\bigr] =
m^{-1} \sum_{\ell=1}^{m} \mathbb E\bigl[{\ln
 p_{\chi}^{\theta}\bigl(Z_{\ell
}|\chunk{Z} {0} {\ell-1}\bigr)}\bigr].
\end{equation}
The Cesaro mean convergence lemma implies that, $\PP{}$-a.s.,
%
\begin{equation}
\label{eqidentity-2} \qquad\limsup_{m \to\infty} m^{-1} \sum _{\ell= 1}^{m}
\ln p_{\chi}^{\theta}\bigl(Z_{0}| \chunk{Z} {-\ell} {-1}\bigr) =
\lim_{\ell\to\infty} \ln p_{\chi}^{\theta} \bigl(Z_{0}|\chunk{Z}
{-\ell} {-1}\bigr) =\ln p^{\theta}\bigl(Z_{0}|\chunk{Z} {-\infty}
{-1}\bigr).
\end{equation}
Combining (\ref{eqfatou}), (\ref{eqliminfPE}), (\ref{eqPEliminf}),
(\ref{eqidentity-1}) and (\ref{eqidentity-2}) yields to
%
\begin{eqnarray}
\label{eqlowBoundLimit} &&\mathbb E\bigl[{ \ln p^{\theta}\bigl(Z_{0}|\chunk{Z}
{-\infty} {-1}\bigr) }\bigr]
\nonumber\\
&&\qquad \geq\limsup_{m \to\infty} m^{-1} \sum
_{\ell=1}^{m} \mathbb E\bigl[{\ln p_{\chi}^{\theta}\bigl(Z_{\ell}|\chunk{Z} {0}
{\ell-1}\bigr)}\bigr]
\\
&&\qquad = \limsup_{m \to\infty} \bigl\{\mathbb E\bigl[{m^{-1}\ln p_{\chi}^{\theta}\bigl(Z_{0}^{m}\bigr)
}\bigr]
- m^{-1} \mathbb E\bigl[{\ln p_{\chi}^{\theta}(Z_0)}\bigr] \bigr\}>-\infty,\nonumber
\end{eqnarray}
where the last bound follows from~\ref{assumlikelihoodDrift}\ref
{itemcondition-minoration} and the minorization
\[
\ln p_{\chi}^{\theta}\bigl(Z_{0}^{m}\bigr) \geq\ln\chi(\Dset) + \sum
_{i=0}^{m} \ln\inf_{x \in\Dset} \Pblock[
\theta]{Z_i}(x,\Dset).
\]
The proof of~\ref{itemexistence-approx-stat-bloc} follows.\vspace*{9pt}

\textit{Proof of}~\ref{itemconv-unif-entropy-stat-bloc}.
According to Proposition \ref
{propforgetting-initial-condition-likelihood} (\ref
{itemforgetting-initial-condition-loglikelihood}), there exists a
random variable $C$ satisfying $\PP[Y]{C<\infty}=1$
such that for all $k \geq K$ and $m\geq0$,
\[
\sup_{\theta\in\Theta}\bigl\llvert\ln p_{\chi}^{\theta}\bigl(Z_k|
\chunk{Z} {-m} {k-1}\bigr) -\ln p_{\chi}^{\theta}\bigl(Z_k|\chunk{Z}
{-m-1} {k-1}\bigr) \bigr\rrvert\leq C \kappa^{k+m},\qquad \PP{}\mbox{-a.s.},
\]
which implies that
\[
\sup_{\theta\in\Theta} \bigl\llvert\ln p_{\chi}^{\theta}\bigl(Z_k|
\chunk{Z} {0} {k-1}\bigr) -\ln p^{\theta}\bigl(Z_k|\chunk{Z} {-
\infty} {k-1}\bigr) \bigr\rrvert\leq C \kappa^{k}/(1-\kappa),\qquad \PP{}
\mbox{-a.s.}
\]
The proof of~\ref{itemconv-unif-entropy-stat-bloc} follows from the
obvious decomposition
%
\begin{eqnarray}\label{eqdecomp-log-stat}
n^{-1}\ln p_{\chi}^{\theta}\bigl(Z_{0}^{n-1}\bigr) &=& n^{-1}
\sum_{k=1}^{n-1} \ln p_{\chi}^{\theta}\bigl(Z_k|\chunk{Z} {0} {k-1}\bigr) + n^{-1} \ln\dens[\chi] [
\theta]{Z_0},
\nonumber\\[-8pt]\\[-8pt]
n^{-1}\ln p^{\theta}\bigl({\chunk{Z} {0} {n-1}|\chunk{Z} {-\infty}
{-1}}\bigr) &=& n^{-1} \sum_{k=0}^{n-1}
\ln p^{\theta}\bigl(Z_k|\chunk{Z} {-\infty} {k-1}\bigr).\nonumber
\end{eqnarray}
The proof of~\ref{itemconv-ponct-entropy-bloc} follows from (\ref
{eqesperance-log-stat-bloc}) and (\ref{eqdecomp-log-stat})
using the Birkhoff theorem; see, for example,~\cite{walters2000},
Theorem 1.14.
\end{pf}
%
\begin{prop}
\label{propbound-conv-unif-entropy-rate}
Assume~\ref{assumlikelihoodDrift}--\ref{assumcontinuity}. Let $\chi$
be a probability measure such that $\chi\in\mathcal{M}(\Dset,r)$
[where $\mathcal{M}(\Dset,r)$ is defined in (\ref{eqmort-de-rire})].

\renewcommand\thelonglist{(\roman{longlist})}
\renewcommand\labellonglist{\thelonglist}
\begin{longlist}
\item\label{itemmajo-unif-dens-stat} For any $\theta_0 \in\Theta$
and any $\rho>0$,
\[
\limsup_{n \to\infty} \sup_{\theta\in\mathcal{B}(\theta_0,\rho)} \frac
{1}{n} \ln p_{\chi}^{\theta}\bigl(Z_{0}^{n-1}\bigr) \leq
\mathbb E\Bigl[{\sup_{\theta\in\mathcal{B}(\theta
_0,\rho)} \ln p^{\theta}\bigl(Z_0|\chunk{Z} {-\infty} {-1}\bigr)}\Bigr],\qquad
\PP{}\mbox{-a.s.}
\]
\item\label{itemsemi-continuity-limit} The function $\theta\mapsto
\PE{\ln\dens[\ ][\theta]{Z_0|\chunk{Z}{-\infty}{-1}}}$ is upper
semi-continuous.
\item\label{itembound-limsup-rho-majo-unif-dens-stat} For any compact
set $\Xi\subset\Theta$, the sequence
$ (\sup_{\theta\in\Xi} \frac{1}{n} \ln\dens[\chi][\theta
]{\chunk{Z}{0}{n-1}} )_{n\geq0}$ converges $\PP{}$-a.s. and
\[
\lim_{n \to\infty} \sup_{\theta\in\Xi} \frac{1}{n} \ln
p_{\chi}^{\theta}\bigl(Z_{0}^{n-1}\bigr) = \sup_{\theta\in\Xi} \mathbb E\bigl[{\ln
p^{\theta}\bigl(Z_0|\chunk{Z} {-\infty} {-1}\bigr)}\bigr],\qquad \PP{}\mbox{-a.s.}
\]
\end{longlist}
\end{prop}
\begin{pf}
\textit{Proof of}~\ref{itemmajo-unif-dens-stat}.
Proposition~\ref{proplimiteRelativeEntropyRateZ}\ref
{itemconv-unif-entropy-stat-bloc} shows that
%
\begin{eqnarray}\label{eqbound-dens-stat}
&&
\limsup_{n \to\infty} \sup_{\theta\in\mathcal{B}(\theta_0,\rho)} \frac
{1}{n} \ln p_{\chi}^{\theta}\bigl(Z_{0}^{n-1}\bigr)
\nonumber\\[-8pt]\\[-8pt]
&&\qquad\leq\limsup_{n \to\infty} \frac{1}{n} \sum_{i=0}^{n-1}
\sup_{\theta\in\mathcal{B}(\theta_0,\rho)} \ln p^{\theta}\bigl(Z_i|\chunk{Z} {-\infty}
{i-1}\bigr),\qquad \PP{}\mbox{-a.s.}
\nonumber
\end{eqnarray}
By (\ref{eqmajo-cond-dens-par-G}), for any $\theta_0 \in\Theta$ and
$\rho>0$,
%
\begin{eqnarray}
\label{eqdouble-bound}
\ln p^{\theta_0}\bigl(Z_0|
\chunk{Z} {-\infty} {-1}\bigr) &\leq&\sup_{\theta\in\mathcal{B}(\theta_0,\rho
)}\ln p^{\theta}\bigl(Z_0|\chunk{Z} {-\infty} {-1}\bigr)
\nonumber\\[-8pt]\\[-8pt]
&\leq&\sum_{i=0}^{r-1} \sup_{\theta\in\Theta}
\ln^+ \bigl|g(\cdot,Y_i)\bigr|_{\infty
},\qquad \PP{}\mbox{-a.s.},
\nonumber
\end{eqnarray}
which shows using (\ref{eqesperance-log-stat-bloc}) and \ref
{assummajo-g} that
\[
\mathbb E\Bigl[{\Bigl\llvert\sup_{\theta\in\mathcal{B}(\theta_0,\rho)}\ln p^{\theta}\bigl(Z_0|
\chunk{Z} {-\infty} {-1}\bigr)\Bigr\rrvert}\Bigr]<\infty.
\]
The Birkhoff theorem therefore implies
%
\begin{eqnarray}\label{eqconv-bound-stat}
&&
\limsup_{n \to\infty} \frac{1}{n} \sum_{i=0}^{n-1}
\sup_{\theta\in\mathcal{B}(\theta_0,\rho)} \ln p^{\theta}\bigl(Z_i|\chunk{Z} {-\infty}
{i-1}\bigr)
\nonumber\\[-8pt]\\[-8pt]
&&\qquad= \mathbb E\Bigl[{\sup_{\theta\in\mathcal{B}(\theta_0,\rho)} \ln p^{\theta}\bigl(Z_0|\chunk{Z} {-
\infty} {-1}\bigr)}\Bigr],\qquad \PP{}\mbox{-a.s.},
\nonumber
\end{eqnarray}
which completes the proof of~\ref{itemmajo-unif-dens-stat}.\vspace*{9pt}

\textit{Proof of}~\ref{itemsemi-continuity-limit}.
First note that
%
\begin{equation}
\label{equpper-bound-lim-sup-esp-lim} \sup_{\theta\in\mathcal{B}(\theta
_0,\rho)} \mathbb E\bigl[{ \ln p^{\theta}\bigl(Z_0|\chunk{Z} {-\infty} {-1}\bigr)}\bigr] \leq
\mathbb E\Bigl[{\sup_{\theta\in\mathcal
{B}(\theta_0,\rho)}
\ln p^{\theta}\bigl(Z_0|\chunk{Z} {-\infty} {-1}\bigr)}\Bigr].
\end{equation}
Now, since under~\ref{assumcontinuity}, for any $m \geq p$, $\PP{}$-a.s.,
the function $\theta\mapsto\ln\dens[\chi][\theta]{Z_0|\chunk
{Z}{-m}{-1}}$ is continuous,
then $\PP{}$-a.s., the function $\theta\mapsto\ln\dens[\ ][\theta
]{Z_0|\chunk{Z}{-\infty}{-1}}$ is continuous as a uniform limit of continuous
functions. Using (\ref{eqdouble-bound}),
\[
\sum_{i=0}^{r-1} \sup_{\theta\in\Theta} \ln^+
\bigl|g(\cdot,Y_i)\bigr|_{\infty} -\sup_{\theta\in\mathcal{B}(\theta_0,\rho)}\ln
p^{\theta}\bigl(Z_0|\chunk{Z} {-\infty} {-1}\bigr)\geq0,
\]
the monotone convergence theorem therefore implies that
%
\begin{eqnarray}\label{eqcont-vois-limit}
\lim_{\rho\downarrow0} \mathbb E\Bigl[{\sup_{\theta\in\mathcal{B}(\theta
_0,\rho
)} \ln p^{\theta}\bigl(Z_0|
\chunk{Z} {-\infty} {-1}\bigr)}\Bigr]
&=& \mathbb E\Bigl[{ \lim_{\rho\downarrow0}
\sup_{\theta\in\mathcal{B}(\theta _0,\rho )} \ln
p^{\theta}\bigl(Z_0|\chunk{Z} {-\infty} {-1}\bigr)}\Bigr]\nonumber\\[-8pt]\\[-8pt]
&=& \mathbb
E\bigl[{\ln p^{\theta_0}\bigl(Z_0|\chunk{Z} {-\infty}
{-1}\bigr)}\bigr]. \nonumber
\end{eqnarray}
Combining (\ref{equpper-bound-lim-sup-esp-lim}) and (\ref
{eqcont-vois-limit}) shows that
\[
\lim_{\rho\downarrow0} \sup_{\theta\in\mathcal{B}(\theta_0,\rho)}
\mathbb E\bigl[{ \ln p^{\theta}\bigl(Z_0|\chunk{Z} {-\infty}
{-1}\bigr)}\bigr] \leq\mathbb E\bigl[{\ln p^{\theta_0}\bigl(Z_0|\chunk{Z} {-\infty}
{-1}\bigr)}\bigr].
\]

\textit{Proof of}~\ref{itembound-limsup-rho-majo-unif-dens-stat}. By
taking the limit of both sides of~\ref{itemmajo-unif-dens-stat} with
respect to $\rho\downarrow0$, (\ref{eqcont-vois-limit})~shows that for
any $\theta_0 \in\Theta$,
%
\begin{equation}
\label{eqlim-boule-tendant-vers-zero}\qquad \lim_{\rho\downarrow0} \limsup_{n
\to\infty} \sup_{\theta\in \mathcal{B}(\theta_0,\rho)} \frac{1}{n} \ln
p_{\chi}^{\theta}\bigl(Z_{0}^{n-1}\bigr) \leq\mathbb E\bigl[{\ln
p^{\theta_0}\bigl(Z_0|\chunk{Z} {- \infty} {-1}\bigr)}\bigr],\qquad \PP{}\mbox{-a.s.}
\end{equation}
Therefore, for any $\delta>0$
and $\theta_0 \in\Xi$, there exists $\rho_{\theta_0}>0$ such that
\[
\limsup_{n \to\infty} \sup_{\theta\in\mathcal{B}(\theta_0,\rho
_{\theta_0})} \frac{1}{n} \ln p_{\chi}^{\theta}\bigl(Z_{0}^{n-1}\bigr)
\leq\mathbb E\bigl[{\ln p^{\theta_0}\bigl(Z_0|
\chunk{Z} {-\infty} {-1}\bigr)}\bigr]+\delta, \qquad\PP{}\mbox{-a.s.}
\]
Since $\Xi$ is compact, by extracting a finite covering, the latter
inequality shows that
\[
\limsup_{n \to\infty} \sup_{\theta\in\Xi} \frac{1}{n} \ln p_{\chi}^{\theta}\bigl(Z_{0}^{n-1}\bigr) \leq\sup_{\theta_0
\in
\Xi}\mathbb E\bigl[{\ln p^{\theta_0}\bigl(Z_0|\chunk{Z} {-\infty} {-1}\bigr)}\bigr]+\delta, \qquad
\PP{}\mbox{-a.s.}
\]
Since $\delta$ is arbitrary, we therefore have
%
\begin{equation}
\label{equpper-bound-lim-compact} \limsup_{n \to\infty} \sup_{\theta\in
\Xi}
\frac{1}{n} \ln p_{\chi}^{\theta}\bigl(Z_{0}^{n-1}\bigr) \leq
\sup_{\theta_0
\in
\Xi}\mathbb E\bigl[{\ln p_{}^{\theta_0}\bigl(Z_0|
\chunk{Z} {-\infty} {-1}\bigr)}\bigr].
\end{equation}
Now, since for any $\theta_0 \in\Xi$,
\[
\sup_{\theta\in\Xi} \frac{1}{n} \ln p_{\chi}^{\theta}\bigl(Z_{0}^{n-1}\bigr)
\geq\frac{1}{n} \ln p_{\chi}^{\theta_0}\bigl(Z_{0}^{n-1}\bigr).
\]
Proposition~\ref{proplimiteRelativeEntropyRateZ}\ref
{itemconv-ponct-entropy-bloc} yields
\[
\liminf_{n \to\infty} \sup_{\theta\in\Xi} \frac{1}{n} \ln
p_{\chi}^{\theta}\bigl(Z_{0}^{n-1}\bigr) \geq\mathbb E\bigl[{\ln
p^{\theta_0}\bigl(Z_0| \chunk{Z} {-\infty} {-1}\bigr)}\bigr],\qquad
\PP{}\mbox{-a.s.}
\]
$\theta_0$ being arbitrary in $\Xi$, we finally obtain
\[
\liminf_{n \to\infty} \sup_{\theta\in\Xi} \frac{1}{n} \ln
p_{\chi}^{\theta}\bigl(Z_{0}^{n-1}\bigr) \geq\sup_{\theta_0 \in\Xi}
\mathbb E\bigl[{\ln p^{\theta_0}\bigl(Z_0|\chunk{Z} {-\infty}
{-1}\bigr)}\bigr], \qquad\PP{} \mbox{-a.s.}
\]
Combining this inequality with (\ref{equpper-bound-lim-compact})
completes the proof.
\end{pf}
%
\begin{theorem} \label{thmconsistanceZ}
Assume~\ref{assumlikelihoodDrift}--\ref{assumcontinuity}. Then, for
any probability measure $\chi\in\mathcal{M}(\Dset,r)$,
\[
\lim_{n \to\infty} d\bigl(\mle{\chi} {nr},\Theta_b^\star
\bigr)=0,\qquad \PP{}\mbox{-a.s.},
\]
where $\Theta_b^\star\subset\Theta$ is defined by $\Theta_b^\star
\triangleq\argmax_{\theta\in\Theta}\PE{\ln\dens[\ ][\theta
]{Z_0|\chunk
{Z}{-\infty}{-1}}}$.
\end{theorem}
\begin{pf}
By Proposition~\ref{propbound-conv-unif-entropy-rate}\ref
{itemsemi-continuity-limit} the function $\theta\mapsto\PE{\ln
\dens
[\ ][\theta]{Z_0|\chunk{Z}{-\infty}{-1}}}$
is upper semi-continuous. Therefore the set $\Theta_b^\star$
is compact as a closed subset of a the compact set $\Theta$ so that for
any $\delta>0$,
$\Xi_\delta=\{\theta\in\Theta; d(\theta,\Theta_b^\star)\geq
\delta\}$
is also a compact set.
In addition, as a upper semi-continuous function, $\theta\mapsto\PE
{\ln\dens[\ ][\theta]{Z_0|\chunk{Z}{-\infty}{-1}}}$ restricted to
$\Xi_\delta$
attains its maximum which implies that
\[
\sup_{\theta\in\Xi_\delta} \mathbb E\bigl[{\ln
p^{\theta}\bigl(Z_0|\chunk{Z} {- \infty}
{-1}\bigr)}\bigr]=\max_{\theta\in \Xi_\delta} \mathbb E \bigl[{\ln
p^{\theta}\bigl(Z_0| \chunk{Z} {-\infty} {-1}\bigr)}\bigr]< \mathbb E\bigl[{\ln
p^{\theta^\star}\bigl(Z_0|\chunk{Z} {-\infty} {-1}\bigr)}\bigr],
\]
where $\theta^\star$ is any point in $\Theta_b^\star$. Combining this
with Proposition~\ref{proplimiteRelativeEntropyRateZ}\ref
{itemconv-ponct-entropy-bloc} yields
\[
\lim_{n \to\infty} \sup_{\theta\in\Xi_\delta} \frac{1}{n} \ln
p_{\chi}^{\theta}\bigl(Z_{0}^{n-1}\bigr) < \mathbb E\bigl[{\ln
p^{\theta^\star}\bigl(Z_0|\chunk{Z} {-\infty} {-1}\bigr)}\bigr],\qquad
\PP{}\mbox{-a.s.}
\]
Using that
\[
\lim_{n \to\infty} \frac{1}{n} \ln p_{\chi}^{\theta^\star}\bigl(Z_{0}^{n-1}\bigr)
 =\mathbb E\bigl[{\ln p^{\theta^\star}\bigl(Z_0|\chunk{Z} {-\infty} {-1}\bigr)}\bigr],\qquad
\PP{}\mbox{-a.s.}
\]
we finally obtain that $\PP{}$-a.s., $\mle{\chi}{n} \in\Xi_\delta$
finitely many times. The proof is
complete.
\end{pf}

\subsection{\texorpdfstring{Proofs of Proposition \protect\ref{proplimiteRelativeEntropyRateY} and Theorem \protect\ref{thmconsistanceY}}
{Proofs of Proposition 1 and Theorem 2}}
We have now all the tools for obtaining the consistency of the MLE as a
byproduct of the results obtained
for the block MLE. We first state and prove the forgetting of the
initial distribution for the predictive filter.
%
\begin{lem} \label{lemforgetting-filter}
Assume~\ref{assumlikelihoodDrift}. Let $0 < \gamma^- < \gamma^+
\leq
1$. Then, for all $\eta> 0$, there exists $\rho_\eta\in(0,1)$ such
that, for all sequence $(z_i)_{i \geq0}$ satisfying
%
\begin{equation}
\label{eqdef-K-2} n^{-1} \sum_{i=0}^{n-1}
\one_{\Kset}(z_i) \geq\max\bigl(1-\gamma^-, \bigl(1+
\gamma^+\bigr)/2 \bigr),
\end{equation}
all $\beta\in\ooint{\gamma^-,\gamma^+}$, all measurable function $f$,
all probability measures $\chi$ and $\chi'$ and all $\theta\in
\Theta$,
\begin{eqnarray*}
&&
\biggl\llvert\frac{\chi\Pblock[\theta]{\chunk{z}{0}{n-1}}f}{\chi
\Pblock
[\theta]{\chunk{z}{0}{n-1}}\one_\Xset}- \frac{\chi'\Pblock
[\theta
]{\chunk{z}{0}{n-1}} f}{\chi'\Pblock[\theta]{\chunk{z}{0}{n-1}}
\one_\Xset} \biggr\rrvert\\
&&\qquad\leq2 \Biggl
\{\rho^{\lfloor n(\beta- \gamma^-) \rfloor} + \frac
{\eta^{\lfloor n (\gamma^+-\beta) \rfloor/2}}{\chi(\Dset)\chi'(\Dset)}
\Biggl[ \prod
_{i=0}^{n-1} D_{z_i}^2 \Biggr]
\Biggr\}|f|_{\infty},
\end{eqnarray*}
where $D_z$ is defined in (\ref{eqdefinition-D-z}).
\end{lem}
\begin{pf}
By Proposition~\ref{propmajo-difference-numerateur},
\begin{eqnarray*}
&&
\biggl\llvert\frac{\chi\Pblock[\theta]{\chunk{z}{0}{n-1}}f}{\chi
\Pblock
[\theta]{\chunk{z}{0}{n-1}} \one_\Xset}- \frac{\chi'\Pblock
[\theta
]{\chunk{z}{0}{n-1}} f}{\chi'\Pblock[\theta]{\chunk{z}{0}{n-1}}
\one_\Xset} \biggr\rrvert\\
&&\qquad=
\frac{|\DDelta[\theta]{\chi,\chi'}{\chunk{z}{0}{n-1}}{f,\one_\Xset}|} {
\chi\Pblock[\theta]{\chunk{z}{0}{n-1}}\one_\Xset\times
\chi'\Pblock[\theta]{\chunk{z}{0}{n-1}}\one_\Xset}
\\
&&\qquad\leq2 \rho^{\lfloor n(\beta- \gamma^-) \rfloor} |f|_{\infty}+2 \eta
^{\lfloor n(\gamma^+-\beta)\rfloor/2}
\frac{\prod_{i=0}^{n-1}
|\Pblock[\theta]{z_i}(\cdot, \Xset)|_{\infty}^2} {
\chi\Pblock[\theta]{\chunk{z}{0}{n-1}}\one_\Xset\times
\chi'\Pblock[\theta]{\chunk{z}{0}{n-1}}\one_\Xset} |f|_{\infty},
\end{eqnarray*}
where we have used that
\[
\frac{\chi\Pblock[\theta]{\chunk{z}{0}{n-1}}f}{\chi\Pblock
[\theta
]{\chunk{z}{0}{n-1}}\one_\Xset} \vee\frac{\chi'\Pblock[\theta
]{\chunk
{z}{0}{n-1}} f}{\chi'\Pblock[\theta]{\chunk{z}{0}{n-1}} \one_\Xset} \leq
|f|_{\infty}.
\]
The proof follows by noting that (\ref{eqmino-denom}) implies that
\[
\frac{\prod_{i=0}^{n-1} |\Pblock[\theta]{z_i}(\cdot, \Xset
)|_{\infty}^2} {
\chi\Pblock[\theta]{\chunk{z}{0}{n-1}}\one_\Xset\times
\chi'
\Pblock[\theta]{\chunk{z}{0}{n-1}}\one_\Xset}\leq\frac{ [
\prod_{i=0}^{n-1} D_{z_i}^2 ]}{\chi(\Dset)\chi'(\Dset)}.\hspace*{75pt}\qed\vadjust{\goodbreak}
\]
\noqed\end{pf}
\begin{pf*}{Proof of Proposition~\ref{proplimiteRelativeEntropyRateY}}
\textit{Proof of}~\ref{itemlim-ponct-cond-rate-entropy-Y}.
Let $\chi$ a probability measure such that $\chi(\Dset)>0$.
The first step of the proof consists of using the forgetting property
obtained in Lemma~\ref{lemforgetting-filter} to show that $\PP{}$-a.s.,
the sequence $(\dens[\chi][\theta]{Y_0|\chunk{Y}{-\ell
}{-1}})_{\ell
\geq0}$ converges.
Denote for any $t \in\{1,\ldots,r\}$,
\[
\chi^{\theta}_{m,t}(\Aset)= \frac{\chi\Pblock[\theta]{\chunk
{y}{-mr-t}{-mr-1}} \one_A}{\chi\Pblock[\theta]{\chunk
{y}{-mr-t}{-mr-1}} \one_\Xset}.
\]
Then, write for any $m\geq0$, $t\in\{1,\ldots,r\}$ and any $\chunk
{y}{-mr-t}{0} \in\Yset^{mr+t+1}$,
\[
 p_{\chi}^{\theta}\bigl(y_0|\chunk{y} {-mr-t} {-1}\bigr)=
p_{\chi^\theta_{m,t}}^{\theta}\bigl(y_0|\chunk{z} {-m} {-1}\bigr)=
\frac{\chi^\theta_{m,t}\Pblock[\theta]{\chunk{z}{-m}{-1}}(\g [\theta
](\cdot,y_0))} { \chi^\theta_{m,t} \Pblock[\theta]{\chunk{z}{-m}{-1}}
(\one_\Xset)}.
\]
Let $0 < \gamma^- < \gamma^+ <1$.
Lemma~\ref{lemforgetting-filter} shows that for any $t \in\{1,\ldots,r\}
$ and $\eta>0$, there exists $\rho\in(0,1)$ such that,
if
\[
m^{-1} \sum_{i=-m}^{-1}
\one_{\Kset}(z_i) \geq\max\bigl(1 - \gamma^-,\bigl(1+
\gamma^+\bigr)/2 \bigr),
\]
then for all $\beta\in\ooint{\gamma^-,\gamma^+}$, and $\theta\in
\Theta$,
\begin{eqnarray*}
&&\bigl| p_{\chi}^{\theta}\bigl(y_0|\chunk{y} {-mr-t} {-1}\bigr)-
p_{\chi}^{\theta}\bigl(y_0|\chunk{y} {-mr} {-1}\bigr)\bigr|
\\
&&\qquad \leq2 \Biggl( \rho^{\lfloor m(\beta-\gamma^-)\rfloor} + \frac{\eta
^{\lfloor
m(\gamma^+-\beta) \rfloor/2}}{\chi^\theta_{m,t}(\Dset) \chi
(\Dset)} \prod
_{j=-m}^{-1} (D_{z_j} )^2
\Biggr)\sup_{\theta\in
\Theta} \bigl|\g[\theta](\cdot,y_0)\bigr|_{\infty}
\\
&&\qquad \leq2 \Biggl( \rho^{\lfloor m(\beta-\gamma^-)\rfloor} + \eta^{\lfloor
m(\gamma^+-\beta) \rfloor/2} D'_{-m}
\prod_{j=-m}^{-1} (D_{z_j}
)^2 \Biggr)\sup_{\theta\in\Theta}\bigl|\g[\theta](\cdot,y_0)\bigr|_{\infty},
\end{eqnarray*}
where
\[
D'_{-m}=\max_{t=1,\ldots,r-1}\frac{1}{\inf_{\theta\in\Theta}\chi^{\theta
}_{m,t}(\Dset) \chi(\Dset)}.
\]
$(D'_{-m})_{m \geq0}$ is a stationary sequence. Using the same
argument as in the proof of~(\ref{eqmajo-ln-D}),
the condition $\chi\in\mathcal{M}(\Dset,r)$ [defined in (\ref
{eqmort-de-rire})], we have $\PE{\ln^+ D'_{-m}}<\infty$.
By choosing $\gamma^+$ and $\gamma^-$ such that $\PP[Y]{Z_{0} \in
\Kset
}>\max(1-\gamma^-,(1+\gamma^+)/2)$ and by applying Lemma \ref
{lemboundFilter}, it follows that there exist $\varrho_{\chi} \in
(0,1)$ and a $\PP{}$-a.s. finite random variable $C_{\chi}$ such that
for any $\ell\geq1$,
\[
\bigl| p_{\chi}^{\theta}\bigl(Y_0|\chunk{Y} {-\ell} {-1}\bigr)-
p_{\chi}^{\theta}\bigl(Y_0|\chunk{Y} {-\ell-1} {-1}\bigr)\bigr|\leq C_\chi
\varrho_\chi^{\ell},\qquad \PP{}\mbox{-a.s.}
\]
Similarly, for any probability measure $\chi'$ such that $\chi'(\Dset
)>0$, there exist $\varrho_{\chi,\chi'} \in(0,1)$ and a
$\PP
{}$-a.s.
finite random variable $C_{\chi,\chi'}$ such that for any $\ell\geq0$,
\[
\bigl| p_{\chi}^{\theta}\bigl(Y_0|\chunk{Y} {-\ell} {-1}\bigr)-
p_{\chi'}^{\theta}\bigl(Y_0|\chunk{Y} {-\ell} {-1}\bigr)\bigr|\leq C_{\chi,\chi'}
\varrho_{\chi,\chi'}^{\ell}, \qquad\PP{}\mbox{-a.s.}
\]
This implies that for any probability measure $\chi$ satisfying
$\chi(\Dset)>0$, the sequence $(\dens[\chi][\theta]{Y_0|\chunk
{Y}{-\ell}{-1}})_{\ell\geq0}$
converges $\PP{}$-a.s. and that the limit denoted by $\dens[\ ][\theta
]{Y_0|\chunk{Y}{-\infty}{-1}}$ does not depend on $\chi$. Then, by
stationarity of $(Y_\ell)_{\ell\in\zset}$, we obtain that for all
$k\geq0$ and $\theta\in\Theta$,
\[
\lim_{m \to\infty} p_{\chi}^{\theta}\bigl(Y_k|\chunk{Y} {-m} {k-1}\bigr)=
p^{\theta}\bigl(Y_k|\chunk{Y} {-\infty} {k-1}\bigr),\qquad \PP{}\mbox{-a.s.},
\]
which shows the first part of \ref
{itemlim-ponct-cond-rate-entropy-Y}. To complete the proof of \ref
{itemlim-ponct-cond-rate-entropy-Y}, it remains to prove that $\PE
{|\ln\dens[\ ][\theta]{Y_k|\chunk{Y}{-\infty}{k-1}}|}<\infty$. Since
$\dens[\chi][\theta]{Y_k|\chunk{Y}{-m}{k-1}}\leq\sup_{x \in\Xset}
g^{\theta}(x,Y_k)$, we have
\[
\ln^+  p_{\chi}^{\theta}\bigl(Y_k|\chunk{Y} {-\infty} {k-1}\bigr) \leq\ln^+
\sup_{x \in\Xset} g^{\theta}(x,Y_k),
\]
which shows, under~\ref{assummajo-g}, that
%
\begin{equation}
\label{eqmajo-esp-ln-plus} \mathbb E\bigl[{\ln^+  p^{\theta}\bigl(Y_k|
\chunk{Y} {-\infty} {k-1}\bigr)}\bigr]<\infty.
\end{equation}
This allows us to define $\PE{\ln\dens[\ ][\theta]{Y_k|\chunk
{Y}{-\infty}{k-1}}}$ as
\[
\mathbb E\bigl[{\ln p^{\theta}\bigl(Y_k|\chunk{Y} {-\infty}
{k-1}\bigr)}\bigr]=\mathbb E\bigl[{\ln^+  p^{\theta}\bigl(Y_k|\chunk{Y}
{-\infty} {k-1}\bigr)}\bigr]-\mathbb
E\bigl[{\ln^-p^{\theta}\bigl(Y_k|\chunk{Y} {-\infty}
{k-1}\bigr)}\bigr],
\]
so that $\mathbb E[{\ln^-  p^{\theta}(Y_k|\chunk{Y}{-\infty
}{k-1})}]<\infty$ provided that we have shown\break $\mathbb E[\ln
p^{\theta}(Y_k|\chunk{Y}{-\infty}{k-1})]>-\infty$. By
stationarity of $(Y_k)_{k\in\zset}$,
%
\begin{eqnarray}\label{eqln-dens-ponct-bloc}\quad
r \mathbb E\bigl[{ \ln p^{\theta}\bigl(Y_0|\chunk{Y} {-\infty}
{-1}\bigr)}\bigr] &=& r \bigl\{\mathbb E\bigl[{ \ln^+  p^{\theta}\bigl(Y_0|\chunk{Y} {-
\infty} {-1}\bigr)}\bigr]-\mathbb E\bigl[{ \ln^-  p^{\theta}\bigl(Y_0|\chunk{Y} {-\infty}
{-1}\bigr)}\bigr] \bigr\}
\nonumber\\
&=& \mathbb E\Biggl[{\sum_{k=0}^{r-1} \ln^+
p^{\theta}\bigl(Y_k|\chunk{Y} {-\infty} {k-1}\bigr)}\Biggr]-\mathbb E\Biggl[{ \sum
_{k=0}^{r-1}\ln^-  p^{\theta}\bigl(Y_k|\chunk{Y} {-\infty} {k-1}\bigr)}\Biggr]
\\
&=& \mathbb E\Biggl[{\sum_{k=0}^{r-1}
\ln p^{\theta}\bigl(Y_k|\chunk{Y} {-\infty} {k-1}\bigr)}\Biggr],\nonumber
\end{eqnarray}
where the last equality follows by applying ${\mathbb E}(A-B)={\mathbb
E}(A)-{\mathbb E}(B)$ for nonnegative random variables $A,B$ such that
${\mathbb E}(A)<\infty$. Now, note that
\begin{eqnarray*}
\prod_{k=0}^{r-1}  p^{\theta}\bigl(Y_k|\chunk{Y} {-\infty}
{k-1}\bigr)&=&\prod _{k=0}^{r-1} \lim_{m \to\infty}
p_{\chi}^{\theta}\bigl(Y_k|\chunk{Y} {-mr} {k-1}\bigr)=\lim_{m \to\infty}\prod
_{k=0}^{r-1}  p_{\chi}^{\theta}\bigl(Y_k| \chunk{Y} {-mr} {k-1}\bigr)
\\
&=&
\lim_{m \to\infty}p_{\chi}^{\theta}\bigl({\chunk{Y} {0} {r-1}|\chunk{Y} {-mr}
{-1}}\bigr)=\lim_{m \to\infty}  p_{\chi}^{\theta}\bigl(Z_0|\chunk{Z} {-m}
{-1}\bigr)\\
&=& p^{\theta}\bigl(Z_0|\chunk{Z} {-\infty} {-1}\bigr).
\end{eqnarray*}
By plugging this expression into (\ref{eqln-dens-ponct-bloc}) and
using $\PE{|\ln\dens[\chi][\theta]{Z_0|\chunk{Z}{-\infty}{-1}}
|}<\infty$ (see Proposition~\ref{proplimiteRelativeEntropyRateZ}), we
finally obtain
%
\begin{equation}
\label{eqegalite-limite-Y-et-limite-Z} r\mathbb E\bigl[{ \ln p^{\theta}\bigl(Y_0|
\chunk{Y} {-\infty} {-1}\bigr)}\bigr]=\mathbb E\bigl[{\ln p^{\theta}\bigl(Z_0|\chunk{Z}
{-\infty} {-1}\bigr)}\bigr] >-\infty,
\end{equation}
which completes the proof of
\ref{itemlim-ponct-cond-rate-entropy-Y}.\vadjust{\goodbreak}

\textit{Proof of}~\ref{itemlim-ponct-rate-entropy-Y}.
Let $\chi$ be a probability measure such that $\chi(\Dset)>0$ and
let $t \in\{0,\ldots,r-1\}$. Then, for any $m\geq0$,
%
\begin{eqnarray}\label{eqdouble-majo-dens-Y-bloc}\quad
m^{-1} \ln p_{\chi}^{\theta}\bigl(Z_{0}^{m+1}\bigr) &\leq&
m^{-1} \ln p_{\chi}^{\theta}\bigl(Y_{0}^{mr+t}\bigr) +
m^{-1} \ln^+ A_{m,t}
\nonumber\\[-8pt]\\[-8pt]
&\leq& m^{-1} \ln p_{\chi}^{\theta}\bigl(Z_{0}^{m}\bigr) +
m^{-1} \ln^+ B_{m,t}+m^{-1} \ln^+
A_{m,t},
\nonumber
\end{eqnarray}
where
\[
A_{m,t}\triangleq\sup_{\theta\in\Theta} \sup_x
p_{Q^\theta(x,\cdot)}^{\theta}\bigl(Y_{mr+t+1}^{(m+1)r-1}\bigr),\qquad
B_{m,t}\triangleq\sup_{\theta\in\Theta} \sup_x  p_{\delta_x}^{\theta}\bigl(Y_{mr}^{mr+t}\bigr).
\]
Note that $(A_{m,t})_{m \geq0}$ and $(B_{m,t})_{m \geq0}$ are
stationary. Moreover, using~\ref{assummajo-g}, it can be easily
checked that
\[
\mathbb E\bigl[{\ln^+ A_{m,t}}\bigr]<\infty,\qquad \mathbb E\bigl[{\ln^+ B_{m,t}}\bigr]<\infty.
\]
Then, Lemma~\ref{lemstationary} may apply and for any $\beta\in
(0,1)$, there exist $\PP{}$-a.s. finite random variables $A,B$ such that
for all $m\geq0$,
\[
A_{m,t}\leq A \beta^{-m},\qquad B_{m,t}\leq B
\beta^{-m},\qquad \PP{}\mbox{-a.s.}
\]
so that, $\PP{}$-a.s.,
\begin{eqnarray*}
&&0 \leq\limsup_{m \to\infty} m^{-1} \ln^+ A_{m,t} \leq-\ln
\beta,
\\
&&0 \leq\limsup_{m \to\infty} m^{-1} \ln^+ B_{m,t} \leq-\ln
\beta.
\end{eqnarray*}
By letting $\beta\uparrow1$,
%
\begin{equation}
\label{eqlim-AnBn} \lim_{m \to\infty}
m^{-1} \ln^+ A_{m,t}=0,\qquad \lim_{m \to\infty} m^{-1}
\ln^+ B_{m,t}=0,\qquad \PP{}\mbox{-a.s.}
\end{equation}
Now, note that $(A_{m,t})_{m \geq0}$ and $(B_{m,t})_{m \geq0}$ do not
depend on $\theta\in\Theta$ so that (\ref{eqlim-AnBn}) together
with (\ref{eqdouble-majo-dens-Y-bloc}) yields
%
\begin{equation}
\label{eqdiff-unif-rate-entropy-Y-Z} \limsup_{m \to\infty}
\sup_{\theta\in\Theta} m^{-1} \bigl|\ln p_{\chi}^{\theta}\bigl(Y_{0}^{mr+t}\bigr)
- \ln p_{\chi}^{\theta}\bigl(Z_{0}^{m}\bigr)\bigr| =0, \qquad\PP{}\mbox{-a.s.}
\end{equation}
Since $t$ is chosen arbitrarily in $\{0,\ldots,r-1\}$, we finally obtain
using Proposition~\ref{proplimiteRelativeEntropyRateZ}\ref
{itemconv-unif-entropy-stat-bloc},
\begin{eqnarray*}
\lim_{n \to\infty} n^{-1} \ln p_{\chi}^{\theta}\bigl(Y_{0}^{n}\bigr)
&=&r^{-1} \lim_{m \to\infty} m^{-1} \ln p_{\chi}^{\theta}\bigl(Z_{0}^{m}\bigr)
\\
&=& r^{-1}\mathbb E\bigl[{\ln p^{\theta}\bigl(Z_0|\chunk{Z} {-\infty}
{-1}\bigr)}\bigr]\\
&=&\mathbb E\bigl[{ \ln p^{\theta}\bigl(Y_0|\chunk{Y} {-\infty}
{-1}\bigr)}\bigr],\qquad
\PP{}\mbox{-a.s.},
\end{eqnarray*}
which completes the proof of
Proposition~\ref{proplimiteRelativeEntropyRateY}.
\end{pf*}
\begin{pf*}{Proof of Theorem~\ref{thmconsistanceY}}
By Proposition~\ref{propbound-conv-unif-entropy-rate}\ref
{itemsemi-continuity-limit} and (\ref
{eqegalite-limite-Y-et-limite-Z}), the function $\theta\mapsto\ell
(\theta)$ is upper semi-continuous. Moreover, (\ref
{eqegalite-limite-Y-et-limite-Z}) also implies
\[
\Theta^\star=\argmax_{\theta\in\Theta}\mathbb E\bigl[{\ln
p^{\theta}\bigl(Y_0|\chunk{Y} {-\infty}
{-1}\bigr)}\bigr]=\argmax_{\theta\in\Theta}\mathbb E\bigl[{\ln
p^{\theta}\bigl(Z_0|\chunk{Z} {-\infty} {-1}\bigr)}\bigr]=\Theta_b^\star.
\]
Now let $t$ in $\{0,\ldots,r-1\}$ and recall that $\chunk
{Z}{0}{m}=\chunk
{Y}{0}{mr-1}$. Theorem~\ref{thmconsistanceZ} together with (\ref
{eqdiff-unif-rate-entropy-Y-Z}) shows that
%
\begin{equation}
\label{eqconsist-nr+t} \lim_{n \to\infty} d\bigl(\mle{\chi} {nr+t},
\Theta^\star\bigr)=0,\qquad \PP{}\mbox{-a.s.}
\end{equation}
The proof of Theorem~\ref{thmconsistanceY} is then complete since $t$
is arbitrary in $\{0,\ldots,r-1\}$.
\end{pf*}
\begin{pf*}{Proof of
Proposition~\ref{propsufficient-condition-mdr}}
Under these two conditions, for any $u \in\{1,\ldots,r\}$, and $\theta
\in\Theta$,
\begin{eqnarray*}
\hspace*{-4pt}&&\chi\mathbf{L}^{\theta}\bigl\langle \chunk{y} {0} {u-1} \bigr\rangle  \one_\Dset
\\
\hspace*{-4pt}&&\quad\geq\Biggl( \prod_{i=0}^{u-1}
\inf_{x_i \in\Dset_i} \g[\theta](x_i,y_i) \Biggr) \int
\cdots\int\chi(\rmd x_0) \one_{\Dset}(x_u) \prod
_{i=1}^{u} \one_{\Dset
_{i-1}}(x_{i-1})
\Q[\theta](x_{i-1},\rmd x_i)
\\
\hspace*{-4pt}&&\quad\geq\Biggl(\prod_{i=0}^{u-1}
\inf_{x_i \in\Dset_i} \g[\theta](x_i,y_i) \Biggr) \chi(
\Dset_0) \delta^u.
\end{eqnarray*}
\upqed\end{pf*}
\begin{pf*}{Proof of Lemma~\ref{lemcheck-CS2-mort-de-rire}}
The proof proceeds by induction on $u \in\{1,\ldots,r\}$. Assume that
$\Dset_{u-1}$ is a compact subset; we show that there exists a compact
set $\Dset_u$ such that $\inf_{x_{u-1} \in\Dset_{u-1}} \inf_{\theta
\in\Theta} \Q[\theta](x_{u-1},\Dset_u) \geq\delta$.

Let $(x,\theta) \in\Dset_{u-1} \times\Theta$ and set $ \delta<
\delta' < 1$.
Since $\Xset= \rset^d$ is a complete separable metric space and
$\Xsigma
$ is the associated Borel $\sigma$-field,
there exists a sequence $B^{x,\theta}_1,B^{x,\theta}_2,\ldots,$ of open
balls of radius $1$ covering $\Xset$. Choose $N_{x,\theta}$ large
enough so that $\Q[\theta](x, \Oset_{x,\theta}) \geq\delta'$, where
$\Oset_{x,\theta} = \bigcup_{i \leq N_{x,\theta}} B^{x,\theta}_i$.
Since for any open set $\Oset$ the function $(x',\theta') \mapsto\Q
[\theta'](x', \Oset)$ is lower semi-continuous, there exists a
neighborhood $\mathcal{V}_{x,\theta}$ (for the product topology on
$\Xset\times\Theta$), such that for all $(x',\theta') \in\mathcal
{V}_{x,\theta}$, $\Q[\theta'](x',\Oset_{x,\theta}) \geq\delta$. Since
$\Oset_{x,\theta}$ is totally bounded
its closure, denoted $\Kset_{x,\theta}$, is a compact subset, which
satisfies, for any $(x',\theta') \in\mathcal{V}_{x,\theta}$
that $\Q[\theta](x,\Kset_{x,\theta}) \geq\delta$.

Then, $\bigcup_{(x,\theta) \in\Dset_{u-1} \times\Theta} \mathcal
{V}_{x,\theta}$ is a covering of $\Dset_{u-1} \times\Theta$. Since the
set $\Dset_{u-1} \times\Theta$ is compact, we may extract a finite
subcover $\Dset_{u-1} \times\Theta\subseteq\bigcup_{i=1}^I
\mathcal
{V}_{x_i,\theta_i}$. Take $\Dset_u= \bigcup_{i=1}^I \Kset_{x_i,\theta
_i}$. As a finite union of compact sets, $\Dset_u$ is a compact set,
which satisfies, for all $(x,\theta) \in\Dset_{u-1} \times\Theta$,
$\Q
[\theta](x,\Dset_u) \geq\delta$. This completes the proof.
\end{pf*}



\printaddresses

\end{document}